\documentclass[aip,superscriptaddress]{elsarticle}
\usepackage[utf8]{inputenc} 
\usepackage[T1]{fontenc}
\usepackage{amsmath,amssymb,amsfonts}
\usepackage{bbold}
\usepackage{color}
\usepackage{url}
\usepackage{graphicx}
\usepackage{tikz}
\usepackage{epstopdf}

\usepackage{hyperref}

\usepackage{amsthm}

\theoremstyle{plain}
\newtheorem{theorem}{Theorem}[section]
\newtheorem{lemma}[theorem]{Lemma}
\newtheorem{corollary}[theorem]{Corollary}
\newtheorem{proposition}[theorem]{Proposition}
\newtheorem{assumptions}[theorem]{Assumptions}

\theoremstyle{definition}

\newtheorem{conjecture}[theorem]{Conjecture}
\newtheorem{remark}[theorem]{Remark}

\newcommand{\diag}{{\rm diag\,}}
\newcommand{\sign}{{\rm sign\,}}

\newcommand{\tr}{{\rm tr\,}}

\newcommand{\U}{{\rm U\,}}

\newcommand{\IM}{{\rm Im\,}}
\newcommand{\eins}{\leavevmode\hbox{\small1\kern-3.8pt\normalsize1}}

\begin{document}

\title{Hard Edge Statistics of Products of P\'olya Ensembles and Shifted GUE's}
\author[1]{Mario Kieburg\corref{cor1}
\fnref{fn1}}\ead{m.kieburg@unimelb.edu.au}
\address[1]{School of Mathematics and Statistics, The University of Melbourne, 813 Swanston Street, Parkville, Melbourne VIC 3010, Australia}

\newcommand{\corr}[1]{{\color{red}#1}}
\newcommand{\tim}[1]{{\color{blue}#1}}


\begin{abstract}
Very recently, we have shown how the harmonic analysis approach can be modified to deal with products of general Hermitian and complex random matrices at finite matrix dimension. In the present work, we consider the particular product of a multiplicative P\'olya ensemble on the complex square matrices and a Gaussian unitary ensemble (GUE) shifted by a constant multiplicative of the identity. The shift shall show that the limiting hard edge statistics of the product matrix is weakly dependent on the local spectral statistics of the GUE, but depends more on the global statistics via its Stieltjes transform (Green function). Under rather mild conditions for the P\'olya ensemble, we prove formulas for the hard edge kernel of the singular value statistics of the P\'olya ensemble alone and the product matrix to highlight their very close similarity. Due to these observations, we even propose a conjecture for the hard edge statistics of a multiplicative P\'olya ensemble on the complex matrices and a polynomial ensemble on the Hermitian matrices.
\end{abstract}

\maketitle

\section{Introduction}\label{sec:intro}

Products and sums of random matrices have seen a revival in the past decade, see the recent reviews~\cite{AI2015,Burda2013}. The reason is two-fold. First and foremost, new applications such as in telecommunications~\cite{AKW2013,KAAC2019,TV2004,WZCT2015}, machine learning~\cite{HN2018,LC2018,TWJTN2019} and quantum information~\cite{CHN2017,Lakshminarayan2013,RSZ2011} require a better comprehension of products of operators. Therefore, it is necessary to understand the generic statistical behaviours of their spectra. Secondly,  new mathematical techniques have been developed and combined with approaches from areas like free probability~\cite{Voiculescu1991,Speicher2011} and harmonic analysis~\cite{Helgason2000,KK2016,KK2019,FKK2017,Kieburg2017,Kieburg2019,KFI2019,KR2016,FZ2019,Faraut2019, Zuber2018}. This provided the rich soil on which new ideas have sprouted to analytically solve more involved structures such as products and sums or even polynomials of random matrices.

Albeit it is widely believed that in the bulk and at the soft edges a product matrix should follow the universal statistics of Gaussian random matrices, indeed this has been proven for the particular case of products of Ginibre matrices where the number of factors stays finite~\cite{LWZ2016} and even when the number grows sublinearly with the matrix size~\cite{ABK2019,LWW2018,ABK2020}, the hard edges show very special behaviour. This has been seen for the singular value statistics of various products of complex squared matrices~\cite{AIK2013,AKW2013,Forrester2014,Liu2017,FIL2018,KKS2016,KS2014,KZ2014} as well as for products including couplings~\cite{Liu2018,AS2018,ACLS2019} and products of real asymmetric matrices~\cite{KFI2019,FILZ2019} and quaternion anti-self-dual matrices~\cite{FILZ2019}. Thence, the hard edge is extremely sensible of what kind of products one studies. Nevertheless, there are basins of attraction for these kernels as it is well-known for the Bessel kernel~\cite{Forresterbook}, though they are incredibly smaller than those in the soft edge statistics and, certainly, tiny compared to the bulk statistics. The reason for this behaviour is that in the bulk the spectrum is extremely stiff so that the eigenvalues are ``squeezed'' between their neighbouring eigenvalues which is extremely strong. In contrast, the hard edge has an additional repulsion from a boundary that strongly influences the closest eigenvalues and may drastically change their statistics. For instance, the number of zero modes of a chiral ensemble appears as a parameter in the Bessel kernel~\cite{Forresterbook}.

Also products of random matrices have a big impact on the hard edge statistics. This influence is born out the fact that the singular value statistics of the single factors play a crucial role and not the eigenvalues. This becomes apparent when considering the harmonic analysis approach~\cite{Helgason2000} for products of matrices exploited for the complex linear group ${\rm Gl}(n,\mathbb{C})$ in~\cite{KK2016,KK2019} and modified to products involving elements of some Lie algebras in~\cite{Kieburg2019,KFI2019}. Singular values have a natural lower boundary at the origin which is the source of the non-trivial effect on the hard edge statistics.

In the present work, we consider the particular product matrix $G(H-n x\eins_n)G^*$, where $x\in\mathbb{R}$ is fixed, $H$ is drawn from a Gaussian unitary ensemble (GUE)~\cite{Forresterbook}, and $G$ is a multiplicative P\'olya ensemble on ${\rm Gl}(n,\mathbb{C})$, see~\cite{FKK2017} and Sec.~\ref{sec:Polya}. Both random matrices are statistically independent. Therefore, our model is related to the one in~\cite{Liu2017} where $G$ has been a product of Ginibre matrices and the shift $x\eins_n$ has been replaced by another constant matrix where most eigenvalues have been symmetrically distributed at $\pm1$ up to a finite number. These slight, but crucial differences change the limiting statistics completely. The matrix $G^*$ is the Hermitian adjoint of $G$ and $n$ is the matrix size. 

The motivation to investigate this model is less born out of an explicit application but results more from the intriguing idea what is going to happen with the hard edge statistics in matrix products when we have on both sides of the edge eigenvalues. It is an exploratory study to get a first feeling for the possible effects and what may impact the local spectral statistics in such constructions. The choice of a shifted GUE for $(H-n x\eins_n)$ and a P\'olya ensemble matrix $G\in{\rm Gl}(n,\mathbb{C})$ is only due to its analytical tractability. Yet, we believe that the rather surprising observations are more general, hence, the Conjecture~\ref{conjecture}. This is the reason why we have introduced the shift $n x\eins_n$ as then we can check what the influence of the local statistics of $(H-n x\eins_n)$ will be on the hard edge statistics of the product matrix. Possible applications, not of the specific model discussed here but of a more general setting with different $G$ and $(H-n x\eins_n)$, can be expected in telecommunications and lattice quantum field theories as products and sums of matrices naturally occur there.

We note that P\'olya ensembles were originally coined polynomial ensembles of derivative type~\cite{KK2016,KK2019,KR2016}; yet, they were re-baptised due to their intimate relation to P\'olya frequency functions~\cite{Polya1913,Polya1915}. The case for $x=0$ and $G$ being a Ginibre matrix or even a product of Ginibre matrices, which is indeed a P\'olya ensemble, too, has been considered in~\cite{FIL2018} and for a shifted GUE matrix $H$ see~\cite{Liu2017}. By means of this model, we want to investigate the impact of the local and macroscopic spectral statistics of the matrix $H$ on the hard edge statistics of the product. For this aim, we make use of the recently derived statistics of the products of P\'olya ensembles on ${\rm Gl}(n,\mathbb{C})$ and polynomial ensembles~\cite{KZ2014} on the Hermitian matrices ${\rm Herm}(n)$ at finite matrix size, see~\cite{Kieburg2019}. Such product matrices satisfy a determinantal point process~\cite{Borodin1998}, which simplifies the analysis a lot since it reduces the whole statistics to a single kernel. We exploit this and reduce our study to the hard edge limit of this kernel for very mild conditions of the P\'olya ensemble.

Our aim is to trace back the hard edge statistics of the product matrix $G(H-n x_0\eins_n)G^*$ in the form of the kernel to the spectral statistics of $G$ and $H$, especially to identify what will have an impact in these statistics. The idea is to use the transformation formula at finite matrix dimension relating the kernel of $H-n x_0\eins_n$ and $G(H-n x_0\eins_n)G^*$. By tilting and deforming the contours we can perform a saddle point analysis. The point is that the deformation has to run around a branch cut which yields one contribution of the new hard edge kernel that can be identified as a transformation of the sine-kernel. A second contribution results from the saddle points and which is a rescaled version of the hard edge kernel of $GG^*$. Surprisingly, only the Green function of the global level density of $H-n x_0\eins_n$ (a shifted Wigner semi-circle) enters in both terms, nothing else of $H$.

The present work is organized as follows. In Sec.~\ref{sec:model}, we recall some common knowledge on the limiting GUE statistics and the less known multiplicative P\'olya ensembles. Therein, we specify the conditions under which the main theorem is true and, additionally, give the hard edge kernel of the squared singular values for the P\'olya matrix $G$, which has been derived for only very specific random matrices and their products before~\cite{AIK2013,AKW2013,Forrester2014,Liu2017,FIL2018,KKS2016,KS2014,KZ2014}. Moreover, we very briefly recall the transformation formula of the kernel from the matrix $H$ to the product matrix $G(H-n x_0\eins_n)G^*$. Here, it is very helpful that the transformation formula is given in terms of a double contour integral. Section~\ref{sec:result} is mainly devoted to the preparation of our main result Theorem~\ref{thm:main}, especially the proper deformations of the contours. In this way, we would like to convey the main ideas and defer the technical details of the proof to Sec.~\ref{sec:proof}, where a case discussion in the value of $x$ is necessary. In Sec.~\ref{sec:conclusio} we summarize our results and state a conjecture when replacing the matrix $H-nx\eins_n$ by a general Hermitian matrix.

As we are interested in the large $n$ limit the matrix dimension is certainly $n>1$ throughout the whole article. This simplifies some subtleties. Moreover, we denote the positive real line (excluding the origin) by $\mathbb{R}_+$ and when we include the origin we write $\mathbb{R}_+^0$.

\section{Preliminaries: Random Matrix Model}\label{sec:model}

We first recall some well-known facts on the GUE, in Subsec.~\ref{sec:GUE}. We especially highlight the fact that it is an additive P\'olya ensemble on ${\rm Herm}(n)$, see~\cite{FKK2017,Kieburg2017}. This allows to employ a particular double contour formula for its kernel that becomes extremely useful. Additive P\'olya ensembles were originally named polynomial ensembles of derivative type~\cite{KR2016}.

In Subsec.~\ref{sec:Polya}, we introduce the P\'olya ensembles~\cite{FKK2017,Kieburg2019} and their properties. In particular, we state the conditions under which Theorem~\ref{thm:main} holds. These conditions also allow us to compute the hard edge statistics of the P\'olya ensemble itself without the product with the GUE. Exactly the same hard edge kernel will appear again in the product matrix $G(H-n x\eins_n)G^*$. The general transformation formula of the finite $n$ kernel from the one of $H-n x\eins_n$ to $G(H-n x\eins_n)G^*$ is briefly recalled in Subsec.~\ref{sec:FiniteN}. This formula has been very recently derived in~\cite{Kieburg2019}.

\subsection{The Gaussian Unitary Ensemble (GUE)}\label{sec:GUE}

The GUE is a random matrix ensemble on the Hermitian matrices ${\rm Herm}(n)$ with a Gaussian as its probability distribution, i.e.,
\begin{equation}\label{probGUE}
	P(H)=2^{-n/2}\left(\frac{n}{\pi}\right)^{n^2/2}\exp\left[-\frac{\tr H^2}{2n}\right]\quad{\rm with}\quad H\in{\rm Herm}(n).
\end{equation}
The joint probability density of the eigenvalues $a=\diag(a_1,\ldots,a_n)$ of this random matrix is~\cite{Forresterbook}
\begin{equation}
	p(a)=\frac{1}{n!}\left(\prod_{j=0}^{n-1}\frac{1}{\sqrt{2\pi n}\, n^j\, j!}\right)\Delta_n^2(a)\exp\left[-\frac{\tr a^2}{2n}\right],
\end{equation}
where $\Delta_n(a)=\prod_{n\geq b>c\geq1}(a_b-a_c)$ is the Vandermonde determinant. We note that the normalization is chosen in such a way that the mean level spacing stays of order one in the bulk of the spectrum when taking $n\to\infty$.

It is well-known~\cite{Forresterbook} that the GUE satisfies a determinantal point process~\cite{Borodin1998}, in particular its $k$-point correlation function
\begin{equation}\label{kpoint}
	R_k(a_1,\ldots,a_k)=\frac{n!}{(n-k)!}\int_{\mathbb{R}^{n-k}}da_{k+1}\cdots da_np(a)=\det[K_n(a_b,a_c)]_{b,c=1,\ldots, k}
\end{equation}
is completely determined by a single kernel $K_n(a_1,a_2)$. For the GUE this kernel has several representations. We choose the one derived in~\cite[Corollary III.3]{Kieburg2017}
\begin{equation}\label{GUEkernel}
	K_n(a_1,a_2)=\oint_{|z'|=1}\frac{dz'}{2\pi i }\int_{-\infty}^{\infty}\frac{dz}{2\pi }\frac{1-(z/z')^n}{z'-z}\exp\left[n\frac{{z'}^2-z^2}{2}+i(a_1z'-a_2 z)\right],
\end{equation}
since it has several advantages which we will exploit when investigating the product matrix $G(H-n x\eins_n)G^*$.
The contour integral of $z'$ runs counter-clockwise around the origin. In~\cite{Kieburg2017}, this expression was derived for more general additive P\'olya ensembles on the Hermitian matrices. The GUE or, later, the shifted GUE matrix $H+n x\eins_n$ with $x\in\mathbb{R}$ fixed is a very particular case of these ensembles.

\begin{remark}
Usually the kernel~\eqref{GUEkernel} is written symmetrically in both entries $a_1$ and $a_2$ and not in a form where it is purely a polynomial in its first entry, as it is here the case. The difference is a factor of $\exp[-(a_1^2-a_2^2)/4n]$ that drops out in the $k$-point correlation function~\eqref{kpoint}. Thence, it is irrelevant for the spectral statistics of $H$. It becomes crucial when considering the product $G(H-n x\eins_n)G^*$ because the transformation formula, discussed in Subsec.~\ref{sec:FiniteN}, takes into account that the kernel is a polynomial of order $n-1$ in its first variable which fixes its ambiguity.
\end{remark}

The concrete case of the level density is related to the $1$-point correlation function as follows
\begin{equation}
	\rho_n(x)=R_1(nx)=K_n(nx,nx).
\end{equation}
The rescaling of $a=nx$ does not only guarantee that $\rho_n$ is normalized but also that  the limiting spectral density, which is the Wigner semi-circle~\cite{Wigner1958}, has a finite support. In the present work, it has the form
\begin{equation}\label{densGUE}
	\rho_{\rm GUE}(x)=\lim_{n\to\infty}\rho_n(x)=\frac{1}{\pi}\sqrt{1-\frac{x^2}{4}},
\end{equation}
meaning its support is the open interval $]-2,2[$. The corresponding Green function is given by
\begin{equation}\label{GreenGUE}
	G(z)=\int_{-2}^2\frac{\rho_{\rm GUE}(x)dx}{z-x}=\frac{z}{2}-\sqrt{\frac{z}{2}-1}\sqrt{\frac{z}{2}+1}\ {\rm for\ all}\ z\in\mathbb{C}\setminus[-2,2]..
\end{equation}
We will encounter this quantity later when studying the hard edge statistics of the product random matrix $G(H-n x\eins_n)G^*$. This can be already seen, when asking for the saddle points of the integrand in Eq.~\eqref{GUEkernel} for $a_1=a_2=nx$ and $n\gg1$, namely those are given by
\begin{equation}\label{saddlepoints}
	z_\pm=\left\{\begin{array}{cl} \displaystyle -i\frac{x\pm i\sqrt{4-x^2}}{2}, & |x|\leq 2, \\ \displaystyle -i\frac{x\pm\sign(x)\sqrt{x^2-4}}{2}, & |x|\geq 2. \end{array}\right.
\end{equation}
Thence, we see that $z_-=-i\lim_{\epsilon\searrow0}G(x+i\epsilon)$ which we will write in short as $z_-=-iG(x)$. In the ensuing discussion on the product, this saddle point still plays a crucial role, especially the modulus $r=|z_-|=|G(x)|$ and its real and imaginary part will appear frequently.

The shift in $H-n x\eins_n$ via the variable $x\in\mathbb{R}$ will be fixed and selects the position of the macroscopic level density where we zoom in. It can be chosen inside the support of the Wigner semi-circle  $\rho$, on its edges or even outside. We will not restrict this. On the contrary, we would like to see what the difference of the hard-edge statistics is when shifting the matrix $H\in{\rm Herm}(n)$ inside the product $G(H-n x\eins_n)G^*$. Surprisingly the whole statistics will be governed by the Green function~\eqref{GreenGUE}.

Finally we would like to point out what is so appealing about the expression~\eqref{GUEkernel}. For this aim, we consider the shifted variables $(a_1,a_2)\to (nx+a_1,nx+a_2)$ that is a zoom into the local statistics at the point $x\in\mathbb{R}$. To evaluate the kernel, we rescale the contours of $z'$ and $z$ by the radius $r=|z_-|=|G(x)|$, i.e., $z'\to r\,z'$ and $z\to r\,z$, and additional shift the contour of $z$ by $i{\rm Im}(z_+)/r$, which is indeed possible as the integrand is entire and has an exponential decay for $|{\rm Re}(z)|\to\infty$. Then, we can evaluate the first term in the factor $[1-(z/z')^n]/[z'-z]$ via the residue theorem and find
\begin{equation}\label{GUE-analysis}
\begin{split}
	K_n(nx+a_1,nx+a_2)=&r\oint_{|z'|=1}\frac{dz'}{2\pi i }\int_{-\infty+i{\rm Im}(z_+)/r}^{\infty+i{\rm Im}(z_+)/r}\frac{dz}{2\pi }\frac{1-(z/z')^n}{z'-z}\\
	&\times\exp\left[nr^2\left(\frac{{z'}^2-z^2}{2}+i\frac{x}{r}[z'-z]\right)+ir\left(a_1z'-a_2 z\right)\right]\\
	=&r\int_{-|{\rm Re}(z_-)|}^{|{\rm Re}(z_-)|}\frac{dz}{2\pi }\exp\left[ir\left(a_1-a_2\right)z\right]\\
	&-r\oint_{|z'|=1}\frac{dz'}{2\pi i }\int_{-\infty+i{\rm Im}(z_+)/r}^{\infty+i{\rm Im}(z_+)/r}\frac{dz}{2\pi }\frac{(z/z')^n}{z'-z}\\
	&\times\exp\left[nr^2\left(\frac{{z'}^2-z^2}{2}+i\frac{x}{r}[z'-z]\right)+ir\left(a_1z'-a_2 z\right)\right].
\end{split}
\end{equation}
The first integral is evidently the famous sine-kernel~\cite{Forresterbook} on the proper domain when $|x|<2$ (Note that then $r=1$) while it vanishes for $|x|\geq2$. The latter term can be shown to vanish as long as $|x|\neq2$. We do not show this here as it is well-known and the computations are standard, it does not differ much from the one in~\cite{Liu2017}. For the edges at $x=\pm2$ this second term yields the contribution of the Airy-kernel~\cite{Forresterbook}  after proper deformation of the $z$ contour. This deformation is necessary since the singularity at $z=z'$ becomes non-integrable when the two contours cross each other at a vanishing angle.

Interestingly, when considering the product $G(H-n x\eins_n)G^*$ at the hard edge, we can essentially perform the same deformations and, rather as a surprise, the edge contribution at $|x|=2$ of the second  term vanishes then, albeit it is crucial for the pure GUE statistics.

\subsection{P\'olya Ensembles on ${\rm Gl}(n)$}\label{sec:Polya}

The second random matrix $G\in{\rm Gl}(n,\mathbb{C})$  involved in the product $G(H-n x\eins_n)G^*$ is a complex invertible matrix that should be distributed by a multiplicative P\'olya ensemble~\cite{FKK2017,Kieburg2019} and is statistically independent of $H$. Such a P\'olya ensemble is first of all a unitary bi-invariant polynomial ensemble on ${\rm Gl}(n,\mathbb{C})$,  in particular its distribution $Q$ is invariant under $Q(G)=Q(U_1GU_2)$ for all $G\in{\rm Gl}(n,\mathbb{C})$ and unitary matrices $U_1,U_2\in\U(n)$. Hence, the eigenvectors are given by the Haar measure on $\U(n)$ and everything specific of the considered random matrix $G$ is encoded in the joint probability distribution of the squared singular values $\lambda=\diag(\lambda_1,\ldots,\lambda_n)\in\mathbb{R}_+^n$, which has the form
\begin{equation}\label{pol.G}
	q(\lambda)=\frac{1}{n!}\Delta_n(\lambda)\frac{\det[w_b(\lambda_c)]_{b,c=1,\ldots,n}}{\det[\mathcal{M}w_b(c)]_{b,c=1,\ldots,n}}>0
\end{equation}
for polynomial ensembles~\cite{KZ2014}.
The weights $\{w_j\}_{j=1,\ldots,n}$ are all $L^1$-functions on $\mathbb{R}_+$ and their moments from order $0$ to order $n-1$ exist. Thus, the Mellin transform
\begin{equation}\label{Mellin}
	\mathcal{M}w_b(s)=\int_0^\infty\frac{d\lambda}{\lambda}w_b(\lambda)\lambda^s
\end{equation}
is well-defined for all complex $s$ with ${\rm Re}(s)\in[1,n]$. A multiplicative P\'olya ensemble, moreover, satisfies the condition
\begin{equation}\label{Polya.G}
	w_b(\lambda_c)=(-\lambda_c\partial_c)^{b-1}\omega^{(n)}(\lambda_c)\quad{\rm for\ all}\quad b=1\ldots,n,
\end{equation}
where we emphasised a possible explicit $n$-dependence by the superscript.
Relation~\eqref{Polya.G} implies for the joint probability density
\begin{equation}\label{Polya.G.b}
	q(\lambda)=\frac{1}{\prod_{j=1}^{n} j!\mathcal{M}\omega(j)}\Delta_n(\lambda)\det[(-\lambda_c\partial_c)^{b-1}\omega(\lambda_c)]_{b,c=1,\ldots,n}.
\end{equation}
Taking the derivatives require that $\omega$ is $(n-1)$-times differentiable, and the positivity of the joint probability density tells us that $\omega\circ\exp(x)=\omega(e^x)$ is a P\'olya frequency function, see~\cite{Polya1913,Polya1915,FKK2017}. Later on, we even require the $n$th derivative of $\omega$ so that we assume that, too. In spite of that, we believe that all of the result do not necessarily need this additional derivative but then the proofs need to be modified and become more technically challenging.

To keep the discussion at a minimum for the product $G(H-n x\eins_n)G^*$, we assume the following.

\begin{assumptions}[Requirements of the P\'olya Ensemble]\label{assumptions}\

We  assume that the analytic continuation of the Mellin transform $\mathcal{M}\omega^{(n)}$ is holomorphic on $\mathbb{C}\setminus(]-\infty,1[\cup]n,\infty[)$
and satisfies the following conditions:
\begin{enumerate}
\item		 for all fixed $s\in\mathbb{C}\setminus]-\infty,1[$ the point-wise limit of the Mellin transform
				\begin{equation}\label{assump0}
				\mathcal{M}\omega^{(\infty)}(s)=\lim_{n>{\rm Re}\,s,\ n\to\infty}\mathcal{M}\omega^{(n)}(s)
				\end{equation}
				 exists, 
\item		there is a constant $\widetilde{C}>0$ so that the function
			\begin{equation}\label{assump1}
				1/\mathcal{M}\omega^{(n)}(s)\leq \widetilde{C}\ {\rm for\ all}\ s\in[1,n]\ {\rm and}\ n\in\mathbb{N},
			\end{equation}
\item	it is dominated  as follows,
			\begin{equation}\label{assump2}
				{\rm sup}_{{\rm arg}(z)=\theta}\ |\mathcal{M}\omega^{(n)}(1+z)|\leq C(\theta)\ {\rm for\ all}\ n\in\mathbb{N},
			\end{equation}
			where $0\leq C(\theta)<\infty$ for all $\theta\in[\pi/2,\pi[$.
\end{enumerate}
\end{assumptions}

\begin{remark}
Considering closely, these conditions seem to be not very restrictive. For all classical ensembles like Ginibre, Jacobi and Cauchy--Lorentz ensembles as well as some of the Muttalib-Borodin ensembles that are P\'olya ensembles (see~\cite{KK2016,KK2019}) a simple rescaling of $\omega^{(n)}(a)\to \omega^{(n)}(a/\xi(n))/\xi(n)$ with a specific scale $\xi(n)>0$ brings them into a satisfying form because the new Mellin transform becomes $[\xi(n)]^{s-1}\mathcal{M}\omega^{(n)}(s)$. These requirements even seem to be intimately related to the fact that the limiting spectrum of $GG^*$ has a hard edge at the origin since the needed scaling $\xi(n)$ seems to be the one of the eigenvalues about the hard edge. In this light it also explains why the inverse Ginibre matrix corresponding to $\omega(a)=a^{-n-1}e^{-1/a}$ cannot be brought into such a form that satisfies any of the two conditions because its Mellin transform $\Gamma(n-s+1)$ always decreases stronger than the exponential scaling $[\xi(n)]^{s-1}$ can compensate and grows super-exponentially for ${\rm Re}\,s\to-\infty$. Albeit these arguments seem to be quite logic, a proof is missing so that one needs to be careful with this observation.
\end{remark}

For the ensuing propositions, theorems and computations we need three functions that are intimately related to the weight $\omega$. First, we define the polynomial
\begin{equation}\label{chi.def}
\chi^{(n)}(z)=\sum_{j=0}^{n-1}\frac{z^j}{\mathcal{M}\omega^{(n)}(j+1)}
\end{equation}
which allows us to define another polynomial
\begin{equation}\label{Iomegan}
J\omega^{(n)}(z')=\oint_{|\tilde{z}|=1}\frac{d\tilde{z}}{2\pi i \tilde{z}}\chi^{(n)}(\tilde{z})\exp\left[- \frac{z'}{\tilde{z}}\right]=\sum_{j=0}^{n-1}\frac{(-z')^j}{j!\mathcal{M}\omega^{(n)}(j+1)}.
\end{equation}
It is one crucial component of the kernel of the determinantal point process corresponding to the singular value statistics of the random matrix $G$.

The condition~\eqref{assump1} directly implies the following lemma which we will exploit later on.

\begin{lemma}\label{lem:pointwiselimJ}\

The point-wise limit of~\eqref{Iomegan},
\begin{equation}\label{Iomegainf}
J\omega^{(\infty)}(z')=\lim_{n\to\infty}J\omega^{(n)}(z')=\sum_{j=0}^{\infty}\frac{(-z')^j}{j!\mathcal{M}\omega^{(\infty)}(j+1)},
\end{equation}
exists for all $z'$ and is entire. Moreover, the function $J\omega^{(n)}$ and its limit $J\omega^{(\infty)}$ are exponentially bounded, in particular exploiting $\widetilde{C}>0$ from~\eqref{assump1} it holds
 \begin{equation}\label{lbound}
  |J\omega^{(n)}(z')|\leq \widetilde{C}e^{|z'|}
 \end{equation}
 for all $z'\in\mathbb{C}$ and $n\in\mathbb{N}$ and, thence, for the limit $n\to\infty$, too.
 
\end{lemma}

\begin{proof}
We first prove the existence of the point-wise limit $J\omega^{(\infty)}$. Indeed, for any point in the complex plane it holds the uniform bound
\begin{equation}\label{Iboundproof}
|J\omega^{(n)}(z)|\leq\sum_{j=0}^{n-1}\frac{|z|^j}{j!\mathcal{M}\omega^{(n)}(j+1)}<\widetilde{C}\sum_{j=0}^{n-1}\frac{|z|^j}{j!}<\widetilde{C}e^{|z|}
\end{equation}
for all $n\in\mathbb{N}$ and $z\in\mathbb{C}$ due to~\eqref{assump1}. Bolzano-Weierstrass implies that the sequence $\{J\omega^{(n)}(z)\}_{n\in\mathbb{N}}$ for any fixed  $z\in\mathbb{C}$ can be split into subsequences for which the limit exists. We will see that all convergent subsequences limit $J\omega^{(\infty)}$ which implies the existence of the point-wise limit.

The bound~\eqref{Iboundproof} already proves the bound~\eqref{lbound}, and it is also the reason why $J\omega^{(n)}$ is an entire function as the radius of convergence is infinitely large which carries over to the point-wise limit $J\omega^{(\infty)}$.

To show the point-wise limit we consider
\begin{equation}
\begin{split}
\left|J\omega^{(\infty)}(z')-J\omega^{(n)}(z')\right|\leq&\sum_{j=0}^{\infty}\frac{|z'|^{j}}{j!}\left|\frac{1}{\mathcal{M}\omega^{(\infty)}(j+1)}-\frac{\chi_{[0,n-1]}(j)}{\mathcal{M}\omega^{(n)}(j+1)}\right|
\end{split}
\end{equation}
with $\chi_{[0,n-1]}(j)$ the indicator function which is only $1$ when $j\in[0,n-1]$ but otherwise vanishes. As the series $\sum_{j=0}^{\infty}|z'|^{j}/j!$ is absolutely convergent for any fixed $z'\in\mathbb{C}$ and the difference in the second term of the sum is bounded by $\widetilde{C}$ we can push the limit $n\to\infty$ into the sum and use the point-wise convergence
\begin{equation}
\lim_{n\to\infty}\frac{\chi_{[0,n-1]}(j)}{\mathcal{M}\omega^{(n)}(j+1)}=\frac{1}{\mathcal{M}\omega^{(\infty)}(j+1)}
\end{equation}
for any $j\in\mathbb{N}_0$ to find
\begin{equation}
\begin{split}
\lim_{n\to\infty}\left|J\omega^{(\infty)}(z')-J\omega^{(n)}(z')\right|=0
\end{split}
\end{equation}
for any fixed $z'\in\mathbb{C}$. This finishes the proof.
 
\end{proof}
 
We will additionally encounter the function
 \begin{equation}\label{Komegan}
 K\omega^{(n)}(z)=\int_0^\infty \frac{d\lambda}{\lambda}\omega^{(n)}(\lambda) \exp\left[-i\frac{z}{\lambda}\right]=\int_{-\infty}^\infty\frac{ds}{2\pi}\mathcal{M}\omega^{(n)}(1+is)\Gamma(1+is)(iz)^{-is-1}
 \end{equation}
 for all $z\in\mathbb{C}$ for which the integral exists. 
 As can be readily checked $ K\omega^{(n)}$ is a positive $L^1$ function on $-i\mathbb{R}_+$ as it is a Mellin convolution of two positive $L^1$ functions.
 
Its point-wise limit is
 \begin{equation}\label{Komegainf2}
 K\omega^{(\infty)}(z)=\lim_{n\to\infty} K\omega^{(n)}(z)=\int_{-\infty}^\infty\frac{ds}{2\pi}\mathcal{M}\omega^{(\infty)}(1+is)\Gamma(1+is)(iz)^{-is-1}
 \end{equation}
on $\{z\in\mathbb{C}|{\rm Im}\,z<0\}$. It directly follows from Lebesgue's dominated convergence theorem and the bound of $|\mathcal{M}\omega^{(n)}(1+is)|$ on $s\in\mathbb{R}$. Note that for this limit we still do not need~\eqref{assump2} but only that $\omega^{(n)}$ is an $L^1$ function because the additional Gamma function and the condition that ${\rm Im}\, z<0$ renders it absolutely integrable. The point-wise limit of $\mathcal{M}\omega^{(n)}(1)$ to $\mathcal{M}\omega^{(\infty)}(1)$ also tells us that $ K\omega^{(\infty)}$ is also a positive $L^1(-i\mathbb{R}_+)$ because its integral is $\mathcal{M}\omega^{(\infty)}(1)<\infty$. Some further properties of~\eqref{Komegainf2}  are derived in the following lemma.

\begin{lemma}\label{lem:Kbound}\

We choose a $\theta\in]0,\pi/2[$. The limit~\eqref{Komegainf2} as well as $ K\omega^{(n)}$ have a holomorphic extension to $\mathbb{C}\setminus(i\mathbb{R}_+^0)$ with the integral representation
 \begin{equation}\label{Komegainf}
 \begin{split}
 K\omega^{(n)}(z)=&\int_{-\infty}^\infty\frac{ds}{2\pi}\mathcal{M}\omega^{(n)}(1+i\exp[i \sign(s)\theta]\, s)\\
 &\times \Gamma(1+i\exp[i \sign(s)\theta]\, s)(iz)^{-i\exp[i \sign(s)\theta]\, s-1},
 \end{split}
 \end{equation}
where the non-holomorphic cut of the root $(iz)^{-i\exp[i \sign(s)\theta]\, s-1}$ is given along $i\mathbb{R}_+^0$. This extension satisfies the exponential bound
 \begin{equation}\label{Kbound}
 \begin{split}
  | K\omega^{(n)}(z)+ K\omega^{(n)}(-\bar{z})|\leq& c\left(\frac{1}{(\big|{\rm ln}|z|\big|+1)|z|}+\exp[\alpha|z|]\right),\\
    | K\omega^{(n)}(z)- K\omega^{(n)}(-\bar{z})|\leq& c\left(\frac{1}{([{\rm ln}|z|]^2+1)|z|}+\exp[\alpha|z|]\right)
  \end{split}
\end{equation}
for any $z\in\mathbb{C}\setminus(i\mathbb{R}_+^0)$ with $c>0$ and $\alpha\geq0$ two constants independent of $n$ so that it also holds for $K\omega^{(\infty)}$.

When ${\rm Im}(z)<0$, we can even set $\theta=0$, and there is a simpler bound of the form $| K\omega^{(n)}(z)|\leq c/[(\pi/2-|\varphi|)|z|]$ where $\varphi={\rm arg}(iz)$.

\end{lemma}

\begin{proof}

The requirement~\eqref{assump2} and the holomorphy of $\mathcal{M}\omega^{(n)}$ on $\{z\in\mathbb{C}| {\rm Re}\,z<1\ {\rm and}\ {\rm Im}\,z\neq0\}$ allow to tilt the contour $s\in\mathbb{R}$ to $\exp[i \sign(s)\theta]\, s$ with $\theta\in]0,\pi/2[$. In this way, the Gamma function $\Gamma[1+i \exp[i \sign(s)\theta]\, s]$ drops off super-exponentially like $\exp[-\chi s\,{\rm ln}(s)]$ with $\chi>0$ for $|s|\to\infty$. 
Once this tilting of the two rays of integrations are done, one can easily see the holomorphy on $\mathbb{C}\setminus(i\mathbb{R}_+^0)$ as the derivative in $z$ can be commuted with the integration due to Lebesgue's dominated convergence theorem (the derivative only yields a polynomial pre-factor in $s$ and the Gamma function still guarantees the integrability) and that the cut of the exponential function $(iz)^{-i\exp[i \sign(s)\theta]\, s-1}$ is along $i\mathbb{R}_+^0$.

For any finite $iz=|z|e^{i\varphi}\in\mathbb{C}\setminus(\mathbb{R}_-^0)$, where $\varphi\in]-\pi,\pi[$, we can bound the function as follows
 \begin{equation}\label{Kbound2}
 \begin{split}
  | K\omega^{(n)}(z)|\leq &\int_{-\infty}^\infty\frac{ds}{2\pi}|\mathcal{M}\omega^{(n)}(1+i \exp[i \sign(s)\theta]\, s)|\\
  &\times |\Gamma(1+i \exp[i \sign(s)\theta]\, s)|\frac{|z|^{\sin(\theta)|s|}e^{\varphi \cos(\theta)s}}{|z|}\\
  \leq& \frac{C(\theta)}{| z|} \int_{-\infty}^\infty\frac{ds}{2\pi} |\Gamma(1+i \exp[i \sign(s)\theta]\, s)|\ |z|^{\sin(\theta)|s|}e^{\varphi \cos(\theta)s}.
  \end{split}
\end{equation}
The latter integral stays finite for any $z\in\mathbb{C}\setminus(i\mathbb{R}_+^0)$ and is $n$-independent. Thence, we need to study its asymptotic behaviour for small $|z|\to0$ and $|z|\to\infty$.

For $|z|\to0$, we consider the two components
 \begin{equation}
 \begin{split}
  | K\omega^{(n)}(z)+ K\omega^{(n)}(-\bar{z})|\leq& \frac{2C(\theta)}{| z|} \int_{-\infty}^\infty\frac{ds}{2\pi} |\Gamma(1+i \exp[i \sign(s)\theta]\, s)|\\
  &\times |z|^{\sin(\theta)|s|}\cosh[\varphi \cos(\theta)s],\\
    | K\omega^{(n)}(z)- K\omega^{(n)}(-\bar{z})|\leq&  \frac{2C(\theta)}{| z|} \int_{-\infty}^\infty\frac{ds}{2\pi} |\Gamma(1+i \exp[i \sign(s)\theta]\, s)|\\
    &\times |z|^{\sin(\theta)|s|}\bigl|\sinh(e^{i{\rm sign}(s)\theta} s\varphi)\bigl|,
  \end{split}
\end{equation}
separately. We rescale $s\to s/{\rm ln}|z|$ in both integrals and take into account that $|z|<1$ for $|z|\approx0$, i.e.,
 \begin{equation}
 \begin{split}
  | K\omega^{(n)}(z)+ K\omega^{(n)}(-\bar{z})|\leq& \frac{2C(\theta)}{\bigl| z\, {\rm ln}|z|\bigl|} \int_{-\infty}^\infty\frac{ds}{2\pi} \left|\Gamma\left(1+i \exp[i \sign(s)\theta]\, \frac{s}{\bigl|{\rm ln}|z|\bigl|}\right)\right|\\
  &\times e^{-\sin(\theta)|s|}\cosh\left[\varphi \cos(\theta)\frac{s}{{\rm ln}|z|}\right],\\
    | K\omega^{(n)}(z)- K\omega^{(n)}(-\bar{z})|\leq&  \frac{2C(\theta)}{\bigl| z\, {\rm ln}|z|\bigl|} \int_{-\infty}^\infty\frac{ds}{2\pi} \left|\Gamma\left(1+i \exp[i \sign(s)\theta]\, \frac{s}{\bigl|{\rm ln}|z|\bigl|}\right)\right|\\
    &\times e^{-\sin(\theta)|s|}\left|\sinh\left[\varphi e^{i{\rm sign}(s)\theta} \frac{s}{\bigl|{\rm ln}|z|\bigl|}\right]\right|.
  \end{split}
\end{equation}
Choosing $|z|$ small enough so that $\bigl|\varphi/{\rm ln}|z|\bigl|<\cot(\theta)$, the integrands are dominated by a constant times $\exp[-\delta |s|]$ with a $\delta>0$ so that we can push the limit $|z|\to0$ into the integral and we obtain
\begin{equation}
\begin{split}
\lim_{|z|\to0} \bigl| z\,{\rm ln}|z|\bigl|\,| K\omega^{(n)}(z)+ K\omega^{(n)}(-\bar{z})|\leq& \left[2C(\theta) \int_{-\infty}^\infty\frac{ds}{2\pi}  e^{-\delta|s|}\right],\\
\lim_{|z|\to0} |z|[{\rm ln}|z|]^2 | K\omega^{(n)}(z)- K\omega^{(n)}(-\bar{z})|\leq&  2C(\theta)|\varphi|\cos(\theta) \int_{-\infty}^\infty\frac{ds}{2\pi}\ e^{-\delta|s|}|s|.
\end{split}
\end{equation}
This explains the first part of the bound in~\eqref{Kbound}.

What remains is to understand its asymptotic behaviour for $|z|\to\infty$. This can be computed for each ray $s>0$ and $s<0$ separately via a saddle point approximation. Since the Gamma function stays non-zero for any finite $s$, $|s|$ has to be large to compete with the exponentially growing term $|z|^{\sin(\theta)|s|}$. Therefore,  we can approximate the Gamma function by Stirling's formula for the saddle point approximation which is
\begin{equation}
\begin{split}
|\Gamma(1+i \exp[i \sign(s)\theta]\, s)|=&\sqrt{2\pi |s|}\biggl|(i \exp[i \sign(s)\theta-1]\, s )^{i \exp[i \sign(s)\theta]\, s}\\
&\times\left(1+\mathcal{O}\left(\frac{1}{1+|s|}\right)\right)\biggl|\\
=&\sqrt{2\pi |s|} \left(\frac{|s|}{e}\right)^{-\sin(\theta)|s|}e^{-(\theta+\pi\sign(s)/2)\cos(\theta)|s|}\\
&\times\left(1+\mathcal{O}\left(\frac{1}{1+|s|}\right)\right).
\end{split}
\end{equation}
The saddle point equation is, then, given by
\begin{equation}\label{saddeqKbound}
\sin(\theta){\rm ln}\left(\frac{|z|}{s_{{\rm sign}(s)}}\right)-\left(\theta+\left[\frac{\pi}{2}-\varphi\right]{\rm sign}(s)\right)\cos(\theta)=0
\end{equation}
which implies the saddle point
\begin{equation}
s_{{\rm sign}(s)}=|z|\exp\left[\left(\left[\varphi-\frac{\pi}{2}\right]{\rm sign}(s)-\theta\right)\cot(\theta)\right]\geq0.
\end{equation}
When we split the integration domain of $s$ into the interval $[-\sqrt{|z|},\sqrt{|z|}]\cup]-\infty,-\sqrt{|z|}[\cup]\sqrt{|z|},\infty[$ and choose $|z|$ big enough the saddle points at $s_\pm$ are the global maximum on the respective intervals $]-\infty,-\sqrt{|z|}[$ and $]\sqrt{|z|},\infty[$ as then the Gamma function is too far away from the poles on the negative real axis so that the the logarithmic derivative of the integrand in $|s|$ becomes monotonous decreasing, see the left hand side of~\eqref{saddeqKbound}. This allows to make the estimate
 \begin{equation}\label{Kbound.p1}
 \begin{split}
  | K\omega^{(n)}(z)| \leq& \frac{C(\theta)}{| z|} \int_{-\sqrt{|z|}}^{\sqrt{|z|}}\frac{ds}{2\pi} |\Gamma(1+i \exp[i \sign(s)\theta]\, s)|\ |z|^{\sin(\theta)|s|}e^{\varphi \cos(\theta)s}\\
  &+\sum_{L=\pm1}\frac{C(\theta)}{| z|} \int_{\sqrt{|z|}}^\infty\frac{ds}{2\pi} |\Gamma(1+i L\exp[i L\theta]\, s)|\ |z|^{\sin(\theta)s}e^{L\varphi \cos(\theta)s}\\
  \leq& \left(C(\theta)\int_{-\infty}^{\infty}\frac{ds}{2\pi} |\Gamma(1+i \exp[i \sign(s)\theta]\, s)|\ e^{\pi \cos(\theta)s}\right) |z|^{\sin(\theta)\sqrt{|z|}-1}\\
  &+\sum_{L=\pm1}\frac{C(\theta)}{\sqrt{| z|}} \int_{1-s_L/\sqrt{|z|}}^\infty\frac{d\delta s}{2\pi} |\Gamma(1+i L\exp[i L\theta](s_L+\sqrt{|z|}\delta s))|\\
  &\times |z|^{\sin(\theta)(s_L+\sqrt{|z|}\delta s)}e^{L\varphi \cos(\theta)(s_L+\sqrt{|z|}\delta s)} .
  \end{split}
\end{equation}
We recall that $\theta\in]0,\pi/2[$.
In the second integral we have substituted $s=s_L+\sqrt{|z|}\delta s$, and now perform the approximation of its integrand
\begin{equation}
\begin{split}
&|\Gamma(1+i L\exp[i L\theta](s_L+\sqrt{|z|}\delta s))|\ |z|^{\sin(\theta)(s_L+\sqrt{|z|}\delta s)}e^{L\varphi \cos(\theta)(s_L+\sqrt{|z|}\delta s)}\\
=&\sqrt{2\pi s_L} \left(\frac{e|z|}{s_L}\right)^{\sin(\theta)s_L}e^{([\varphi-\pi/2)L-\theta)\cos(\theta)s_L-\sin(\theta)|z|\delta s^2/(2s_L)}\left(1+\mathcal{O}\left(\frac{1}{1+\sqrt{|z|}}\right)\right).
\end{split}
\end{equation}
As $s_L\propto|z|\gg\sqrt{|z|}$  for $|z|\to\infty$, we can extend the integral to the whole real axis which leads to an exponentially small correction and carry out the Gaussian integral, i.e.
 \begin{equation}\label{Kbound.p2}
 \begin{split}
  | K\omega^{(n)}(z)|
  \leq& c_0(\theta) |z|^{\sin(\theta)\sqrt{|z|}-1}+\sum_{L=\pm1}c_L(\theta)e^{\alpha_L(\theta)|z|}\left(1+\mathcal{O}\left(\frac{1}{1+\sqrt{|z|}}\right)\right),
  \end{split}
\end{equation}
where the constants are
\begin{equation}
\begin{split}
c_0(\theta)=&C(\theta)\int_{-\infty}^{\infty}\frac{ds}{2\pi} |\Gamma(1+i \exp[i \sign(s)\theta]\, s)|\ e^{\pi \cos(\theta)s},\\
c_L(\theta)=&\sqrt{\frac{1}{\sin(\theta)}}C(\theta)\exp\left[\left(\frac{(2-L)\pi}{2}-\theta\right)\cot(\theta)\right],\\
\alpha_L(\theta)=&\sin(\theta)\exp\left[\left(\frac{(2-L)\pi}{2}-\theta\right)\cot(\theta)\right]
\end{split}
\end{equation}
as can be seen $L=-1$ grows the fastest of all three terms yielding the second term in the bound~\eqref{Kbound}.

When ${\rm Im}(z)<0$, the integrand~\eqref{Komegan} is indeed absolutely integrable for $\theta=0$ because $\varphi\in]-\pi/2,\pi/2[$ and the integrand is bounded by a constant times the exponential function $\sqrt{|s|+1}\exp[-(\pi/2-{\rm sign}(s)\varphi) |s|]$ as $|\Gamma(1+is)|<\sqrt{2\pi(|s|+1)}e^{-\pi|s|/2}$. The bound of the function follows then
\begin{equation}
\begin{split}
 | K\omega^{(n)}(z)| \leq& \frac{C(\theta)}{| z|}\int_{-\infty}^\infty\frac{ds}{\sqrt{2\pi}} |\Gamma(1+i \, s)|\ e^{\varphi s}\\
 \leq& \frac{C(\theta)}{| z|}\int_{-\infty}^\infty ds\sqrt{|s|+1}|\exp[-(\pi/2-{\rm sign}(s)\varphi) |s|]\\
 \leq&\frac{c}{(\pi/2-|\varphi|)|z|}
 \end{split}
\end{equation}
with $c>0$ some constant. This finishes the second claim .

\end{proof}

On the half line $i\mathbb{R}_+$ the convergence of the integral~\eqref{Komegan} is also given. Yet, we have to respect the cut of the map $z\mapsto(iz)^{-is}$ which might lead to two different values depending on whether one approaches the imaginary axes from the left or the right. Indeed, what we need is the limit of the difference
\begin{equation}\label{Jomegan}
\tilde{J}\omega^{(n)}(y)=\lim_{\epsilon\to0}\frac{i}{2\pi}[K\omega^{(n)}(iy+\epsilon)-K\omega^{(n)}(iy-\epsilon)]
\end{equation}
for $y>0$. For $y<0$, the limit vanishes.
It can be evaluated by employing its integral representation~\eqref{Komegainf2} combined with the tilt of the contour, i.e.,
\begin{equation}
\begin{split}
\tilde{J}\omega^{(n)}(y)=&\lim_{\epsilon\to0}\int_{-\infty}^\infty\frac{ds}{2\pi}\mathcal{M}\omega^{(n)}(1+i \exp[i \sign(s)\theta]\, s)\Gamma(1+i\exp[i \sign(s)\theta]\, s)\\
&\times\frac{i}{2\pi}\left[(-y+i\epsilon)^{-i\exp[i \sign(s)\theta]\, s-1}-(-y-i\epsilon)^{-i\exp[i \sign(s)\theta]\, s-1}\right]\\
=&\lim_{\epsilon\to0}\int_{-\infty}^\infty\frac{ds}{2\pi}\mathcal{M}\omega^{(n)}(1+i\exp[i \sign(s)\theta]\, s)\Gamma(1+i\exp[i \sign(s)\theta]\, s)\\
&\times\frac{\sin\left[(\pi -\beta_\epsilon)(1+i\exp[i \sign(s)\theta]\, s)\right]}{\pi}(y^2+\epsilon^2)^{-(1+i\exp[i \sign(s)\theta]\, s)/2}
\end{split}
\end{equation}
with $-y\pm i\epsilon=\sqrt{y^2+\epsilon^2}\exp[\pm i(\pi -\beta_\epsilon)]$, meaning $\beta_\epsilon=\arctan[\epsilon/y]$. Due to the absolute convergence that has been guaranteed by the finite angle $\theta\in]0,\pi/2[$, we can shift the limit into the integral and get the factor $\sin\left[\pi(1+i\exp[i \sign(s)\theta]\, s)\right]/\pi$. This term can be rewritten into two Gamma functions with the aid of Euler's reflection formula $\Gamma[1-z]\Gamma[z]=\pi/\sin(\pi z)$. One of these Gamma functions cancels with the one in the numerator and the other leads to the result
\begin{equation}\label{Jomegan.b}
\tilde{J}\omega^{(n)}(y)=\int_{-\infty}^\infty\frac{ds}{2\pi}\frac{\mathcal{M}\omega^{(n)}(1+i\exp[i \sign(s)\theta]\, s)}{\Gamma(-i\exp[i \sign(s)\theta]\, s)}y^{-1-i\exp[i \sign(s)\theta]\, s}.
\end{equation}
Let us emphasize that the tilt of the two rays $s>0$ and $s<0$ is still crucial otherwise it might happen that the integral does not exist when $\mathcal{M}\omega^{(n)}$ does not drop off fast enough. Additionally, it allows us to pull the limit $n\to\infty$ into the integral so that we have
\begin{equation}\label{Jomegan.c}
\tilde{J}\omega^{(\infty)}(y)=\lim_{n\to\infty}\tilde{J}\omega^{(n)}(y)=\int_{-\infty}^\infty\frac{ds}{2\pi}\frac{\mathcal{M}\omega^{(\infty)}(1+i\exp[i \sign(s)\theta]\, s)}{\Gamma(-i\exp[i \sign(s)\theta]\, s)}y^{-1-i\exp[i \sign(s)\theta]\, s}.
\end{equation}
Also  $\tilde{J}\omega^{(n)}$ and $\tilde{J}\omega^{(\infty)}$ satisfy a bound which is a direct consequence of Lemma~\ref{lem:Kbound}.

\begin{corollary}\label{cor:Jtildebound}\

The functions $\tilde{J}\omega^{(n)}$ and $\tilde{J}\omega^{(\infty)}$ can be bounded on any bounded interval $]a,b]\subset\mathbb{R}_+$ like
\begin{equation}\label{Jtildebound}
|\tilde{J}\omega^{(n)}(y)|\leq c\frac{1}{([{\rm ln}(y)]^2+1)y}
\end{equation}
with $c$ an $n$-independent constant only depending on the half-open interval $]a,b]$ so that this inequality also holds for $\tilde{J}\omega^{(\infty)}$. Hence, the two functions are absolutely integrable on each bounded interval.
\end{corollary}

Especially the absolute integrability on a bounded interval that has $0$ as its lower boundary will be of use for us later on.

\begin{proof}

We exploit the second equation in~\eqref{Kbound} which holds for any $z=i y+\epsilon$ with $y,\epsilon>0$ regardless how small $\epsilon$ is. Therefore, we can take the limit $\epsilon\to0$ in this inequality leading to the claim.
\end{proof}

\begin{remark}

One last word on  the notation  $J\omega^{(\infty)}$,  $\tilde{J}\omega^{(\infty)}$ and  $K\omega^{(\infty)}$ before going on. It is reminiscent to the case when choosing $G$ as a complex Ginibre matrix, meaning $\omega(a)=e^{-a}$. Then, $J\omega^{(\infty)}$ and $\tilde{J}\omega^{(\infty)}$ become essentially the Bessel functions of the first kind and $K\omega^{(\infty)}$ stands for the modified Bessel functions of the second kind. Those are usually denoted by $J_\nu$ and $K_\nu$. In particular, the relation~\eqref{Jomegan} between $J_\nu$ and $K_\nu$ is well known, e.g., see~\cite[Chapter 9]{ASbook}.
\end{remark}

The functions $J\omega^{(\infty)}$ and  $\tilde{J}\omega^{(\infty)}$ are also encountered in the hard edge statistics of the singular values of the matrix $G$. In~\cite{KK2016,KK2019} we have shown that the kernel at finite $n$ for the eigenvalues of $GG^*$ is equal to
\begin{equation}\label{kernelG}
K_n^{(G)}(\lambda_1,\lambda_2)=\int_0^1dt p_{n-1}(\lambda_1 t)q_n(\lambda_2 t)
\end{equation}
with the polynomial
\begin{equation}\label{polynomial.new}
p_{n-1}(\lambda)=\sum_{j=0}^{n-1}\frac{(n-1)!(-\lambda)^j}{j!(n-1-j)!\mathcal{M}\omega^{(n)}(j+1)}=(n-1)!\oint_{|z|=1}\frac{dz}{2\pi i z^{n+1}}J\omega^{(n)}(z\lambda)e^{z},
\end{equation}
where the contour around the origin is counter-clockwise, and the weight
\begin{equation}
q_n(\lambda)=\lim_{\epsilon\to0}\int_{-\infty}^\infty \frac{ds}{2\pi} (-1)^n\frac{\Gamma[1+is]\mathcal{M}\omega^{(n)}(1+is)}{(n-1)!\Gamma[1+is-n]}\lambda^{-1-is}e^{-\epsilon s^2}.
\end{equation}
The limit $\epsilon\to0$ for $q_n$ is a regularisation which is needed when $\mathcal{M}\omega^{(n)}$ does not drop off fast enough. Note that we employ different normalizations for the weights and polynomials compared to those in~\cite{KK2016,KK2019} in foresight of the limit $n\to\infty$, and also the regularisation is chosen to be a Gaussian which is more convenient when deforming the contour. The whole kernel is still the same. The expression for the weights $q_n$ can be still massaged by employing Euler's reflection formula for both Gamma functions so that the expression simplifies to
\begin{equation}\label{qint}
\begin{split}
q_n(\lambda)=&\lim_{\epsilon\to0}\int_{-\infty}^\infty \frac{ds}{2\pi} (-1)^n\frac{\sin[\pi(1+is-n)]\Gamma[n-is]\mathcal{M}\omega^{(n)}(1+is)}{(n-1)!\sin[\pi(1+is)]\Gamma[-is]}\lambda^{-1-is}e^{-\epsilon s^2}\\
=&\lim_{\epsilon\to0}\int_{-\infty}^\infty \frac{ds}{2\pi} \frac{\Gamma[n-is]\mathcal{M}\omega^{(n)}(1+is)}{\Gamma[n]\Gamma[-is]}\lambda^{-1-is}e^{-\epsilon s^2}.
\end{split}
\end{equation}
When restricting to $\lambda\in]0,1[$, we can even omit the Gaussian regularization at the expense of tilting the contour again into two rays $s\to\exp[i \sign(s)\theta]s$ with $\theta\in]0,\pi/2[$, namely the integrand~\eqref{qint} can be bounded by a constant times $(|s|^n+1)\exp[{\rm ln}(\lambda)\sin(\theta)|s|]$ which is integrable for $\lambda\in]0,1[$. The weight function reads in this way
\begin{equation}\label{weight.new}
\begin{split}
q_n(\lambda)=&\int_{-\infty}^\infty \frac{ds}{2\pi} \frac{\Gamma[n-i\exp[i \sign(s)\theta]s]\mathcal{M}\omega^{(n)}(1+i\exp[i \sign(s)\theta]s)}{\Gamma[n]\Gamma[-i\exp[i \sign(s)\theta]s]}\lambda^{-1-i\exp[i \sign(s)\theta]s}\\
=&\prod_{l=0}^{n-1}\left(1+\frac{1}{l}\partial_\lambda\lambda\right)\omega^{(n)}(\lambda).
\end{split}
\end{equation}
The last line has been already derived in~\cite[Lemma 4.2]{KK2016}. It follows from the definition of the inverse Mellin transform, e.g., see~\cite[Lemma 2.6]{KK2016}, and the relation $\mathcal{M}^{-1}[-isf(s)](\lambda)=\partial_\lambda\lambda \mathcal{M}^{-1}[f(s)](\lambda)$ for a Mellin transform $f$ of an $L^1$-function. It  nicely shows that the moments of $\lambda^jq_n(\lambda)$ exist for all $j=0,\ldots, n-1$ and that its moments are $\Gamma[n-j]\mathcal{M}\omega^{(n)}(1+j)/(\Gamma[n]\Gamma[j])$.

Equations~\eqref{polynomial.new} and~\eqref{weight.new} can be plugged into the integral~\eqref{kernelG} so that the integral over $t$ can be carried out,
\begin{equation}\label{kernelG.b}
\begin{split}
&K_n^{(G)}(\lambda_1,\lambda_2)
=\sum_{j=0}^{n-1}\frac{(n-1)!}{j!(n-1-j)!\mathcal{M}\omega^{(n)}(j+1)}(-\lambda_1)^j\\
\times&\int_{-\infty}^\infty \frac{ds}{2\pi} \frac{\Gamma[n-i\exp[i \sign(s)\theta]s]\mathcal{M}\omega^{(n)}(1+i\exp[i \sign(s)\theta]s)}{\Gamma[n]\Gamma[-i\exp[i \sign(s)\theta]s]}\frac{\lambda_2^{-1-i\exp[i \sign(s)\theta]s}}{j-i\exp[i \sign(s)\theta]s}.
\end{split}
\end{equation}
The pole at $j-i\exp[i \sign(s)\theta]s=0$ is only apparent and cancels with a zero of $1/\Gamma[-i\exp[i \sign(s)\theta]s]$.

Let us emphasize that the integrability in $s$ is guaranteed by the factor $\lambda_2^{-1-i\exp[i \sign(s)\theta]s}$ and not by $\mathcal{M}\omega^{(n)}$ for all $n$. This is extremely important for the proof in the following proposition that tells us what the hard edge limit of the considered P\'olya ensembles is.

\begin{proposition}[Hard Edge Kernel of P\'olya Ensembles]\label{prop:HEKPolya}\

The hard edge limit of the kernel~\eqref{kernelG} for P\'olya ensembles satisfying the conditions~\eqref{assump1} and~\eqref{assump2} is given by the point-wise limit
\begin{equation}\label{hardkernPolya}
K_\infty^{(G)}(y_1,y_2)=\lim_{n\to\infty}\frac{1}{n}K_n^{(G)}\left(\frac{y_1}{n},\frac{y_2}{n}\right)=\int_0^1dt J\omega^{(\infty)}(y_1 t)\tilde{J}\omega^{(\infty)}(y_2 t)
\end{equation}
for any fixed $y_1\in\mathbb{R}_+^0$ and $y_2\in\mathbb{R}_+$.
\end{proposition}

This proposition has been proven for several particular cases of products of complex random matrices~\cite{AIK2013,AKW2013,Forrester2014,Liu2017,FIL2018,KKS2016,KS2014,KZ2014}. The kernel will also reappear in the hard edge statistics of the product $G(H-nx\eins_n)G^*$. The inversion of this statement, namely when the limit exists then we get the two assumptions~\eqref{assump1} and~\eqref{assump2}, would correspond to the aforementioned observation that they equivalently correspond to the existence of a hard edge limit at the origin. This is, however, an open problem yet.

\begin{remark}\label{remark:toprophardegePolya}

As a side-remark, Proposition~\eqref{prop:HEKPolya} also shows that the function $\tilde{J}\omega^{(\infty)}(t)$ is integrable on any finite interval on $\mathbb{R}_+^0$ as we can set $y_1=0$ and choose $y_2>0$ arbitrarily.
\end{remark}

\begin{proof}
We start from the last line in~\eqref{kernelG.b} with $(\lambda_1,\lambda_2)=(y_1,y_2)/n$ and would like to show that Lebesgue's dominated convergence theorem applies.
The $n$-dependent prefactor in the sum,
\begin{equation}
\frac{(n-1)!}{(n-1-j)! n^j}=\prod_{l=1}^{j}\left(1-\frac{l}{n}\right)\leq1,
\end{equation}
is evidently bounded. For the ratio of the Gamma functions in the integral we employ Stirling's formula,
\begin{equation}
\begin{split}
&\left|\frac{\Gamma[n-i\exp[i \sign(s)\theta]s]n^{i\exp[i \sign(s)\theta]s}}{\Gamma[n]}\right|\\
=&\sqrt{\frac{|n-i\exp[i \sign(s)\theta]s|}{n}}\left|\left(1-i\exp[i \sign(s)\theta]\frac{s}{n}\right)^{n-i\exp[i \sign(s)\theta]s}\right|e^{-\sin(\theta)|s|}\\
&\times\exp\left[\mathcal{O}\left(\frac{1}{n}+\frac{1}{|n-i\exp[i \sign(s)\theta]s|}\right)\right].
\end{split}
\end{equation}
The error is uniform since $|n-i\exp[i \sign(s)\theta]s|>0$. Thus, it  is
\begin{equation}
\begin{split}
&\left|\frac{\Gamma[n-i\exp[i \sign(s)\theta]s]n^{i\exp[i \sign(s)\theta]s}}{\Gamma[n]}\right|\\
\leq&C\sqrt{1+|s|}\left(1+\frac{2\sin(\theta)|s|}{n}+\frac{s^2}{n^2}\right)^{(n+\sin(\theta)|s|)/2}\\
&\times\exp\left[-\cos(\theta){\rm arctan}\left(\frac{\cos(\theta)s}{n+\sin(\theta)|s|}\right)s-\sin(\theta)|s|\right]
\end{split}
\end{equation}
for some constant $C$. The ${\rm arctan}$ term behaves approximately linear for large $|s|$ and stays finite for all $s\in\mathbb{R}$. The same also holds for the term
\begin{equation}
n\,{\rm ln}(\sqrt{1+2\sin(\theta)|s|/n+s^2/n^2})<|s|.
\end{equation}
Both are subleading and do not influence the absolute convergence given by $1/|\Gamma[1-i\exp[i \sign(s)\theta]s]|$, where we shifted the argument to avoid the zero. Its exponential leading behaviour goes with the rate $-\sin(\theta)|s|\,{\rm ln}(\sqrt{1+2\sin(\theta)|s|+s^2})$.

 What is really competing with this  leading asymptotic behaviour of $1/|\Gamma[1-i\exp[i \sign(s)\theta]s]|$ is the term $\exp[\sin(\theta)|s|{\rm ln}(\sqrt{1+2\sin(\theta)|s|/n+s^2/n^2})]$. Looking for the supremum of the function
\begin{equation}
f\left(|s|\right)=\frac{1+2\sin(\theta)|s|/n+s^2/n^2}{1+2\sin(\theta)|s|+s^2},
\end{equation}
we see that the quotient of both exponentials is bounded by the identity because the numerator is always smaller than the denominator. Moreover, it is a strictly decreasing function in $n$. For large $|s|$, the quotient
\begin{equation}
\exp[\sin(\theta)|s|{\rm ln}(\sqrt{1+2\sin(\theta)|s|/n+s^2/n^2})]/|\Gamma[1-i\exp[i \sign(s)\theta]s]|
\end{equation}
behaves like $\exp[-{\rm ln}(n)\sin(\theta)|s|]$.
Therefore, we can choose an $n_0>0$ so that for all $n>n_0$ we have
\begin{equation}
\frac{\exp[\sin(\theta)|s|{\rm ln}(\sqrt{1+2\sin(\theta)|s|/n+s^2/n^2})]}{|\Gamma[1-i\exp[i \sign(s)\theta]s]|}\leq \gamma\exp[-{\rm ln}(n_0)\sin(\theta)|s|]
\end{equation}
with some $n$-independent constant $\gamma>0$.
 As the leading contributions of the other terms, including $|y_2^{-i\exp[i \sign(s)\theta]s}|$, asymptote to $\exp[\gamma_1s+\gamma_2]$ with two $n$-independent constants $\gamma_1,\gamma_2\in\mathbb{R}$ (note that $\mathcal{M}\omega^{(n)}$ is bounded by a constant, see~\eqref{assump2}), we can choose $n$ large enough to dominate the integration. In particular we choose an $n_0>0$ and $n>n_0$, so that $\sin(\theta){\rm ln}(n_0)$ is larger than $\gamma_1$. Since the whole integrand is continuous for all $n$ and remains continuous on the integration domain for all $n$ as we have seen in the above discussion, we can replace it by $\tilde{\gamma}_1\exp[-\sin(\theta)({\rm ln}(n_0)-\tilde{\gamma}_2)|s|]$ with $\tilde{\gamma}_1>0$ and $\tilde{\gamma}_2$ two constants and $n>n_0>e^{\tilde{\gamma}_2}$, showing its absolute integrability.

The point-wise limit can be done by Stirling's formula
\begin{equation}
\frac{\Gamma[n-i\exp[i \sign(s)\theta]s]}{\Gamma[n]}=n^{-i\exp[i \sign(s)\theta]s}\left[1+O(n^{-1})\right].
\end{equation}
 After rewriting $1/(j-\exp[i \sign(s)\theta]s)=\int_0^1dt\, t^{j-i\exp[i \sign(s)\theta]s-1}$ in~\eqref{kernelG.b}, we arrive at the claim~\eqref{hardkernPolya}.
\end{proof}

\subsection{Statistics of the Random Matrix Product at Finite $N$}\label{sec:FiniteN}

In~\cite{Kieburg2019}, we have proven that any product random matrix  $GXG^*$  with $X$ a polynomial ensemble on the Hermitian matrices and $G$ a P\'olya ensemble on ${\rm Gl}(n,\mathbb{C})$ yields again a polynomial ensemble on the Hermitian matrices. Say $K_n$ has been the kernel for the matrix $X$, then the new kernel has the simple form~\cite[Proposition IV.8]{Kieburg2019}
\begin{equation}\label{kerneltrafo}
\tilde{K}_{n}(\tilde{a}_1,\tilde{a}_2)=\oint \frac{d\tilde{z}}{2\pi i \tilde{z}}\chi^{(n)}(\tilde{z})\int_{0}^\infty  \frac{ da}{a}\omega^{(n)}\left(a\right)K_{n}\left(\frac{\tilde{a}_1}{\tilde{z}},\frac{\tilde{a}_2}{a}\right).
\end{equation}
Please, recall the definition~\eqref{chi.def} of $\chi^{(n)}$.
The contour of $\tilde{z}$ encircles the origin counter-clockwise. This formula is certainly also true for the case of $X=H-nx\eins_n$ with $H$ being a GUE and $x\in\mathbb{R}$.
Similar transformation formulas have been derived for sums of random matrices~\cite{Kieburg2017,
CKW2015,Kuijlaars2016,KR2016} as well as otherkinds of  products of random matrices~\cite{KK2019,CKW2015,Kuijlaars2016,KKS2016}.

When applying Eq.~\eqref{kerneltrafo} onto the kernel~\eqref{GUEkernel} with $(a_1,a_2)\to(a_1,a_2)+nx\eins_n$, we cannot easily interchange the order of the integrals in any representation of the contour integral. For this we need to start with the original contour for $z$ and shift it slightly to $z\in\mathbb{R}-i{\rm sign}(\tilde{a}_2)\epsilon$ with $\epsilon>0$. This shift guarantees that the exponential function $\exp[-i\tilde{a}_2z/a]$ becomes absolutely integrable in $a$ about the origin $a\to0$. The integrability for $a\to\infty$ is given by $\omega^{(n)}$. Applying the definition~\eqref{Iomegan} and~\eqref{Komegan}, the kernel for $G(H-nx\eins_n)G^*$ is
\begin{equation}\label{kernelprod}
\begin{split}
\tilde{K}_{n}(\tilde{a}_1,\tilde{a}_2)=&\oint_{|z'|=1}\frac{dz'}{2\pi i }\int_{-\infty-i{\rm sign}(\tilde{a}_2)\epsilon}^{\infty-i{\rm sign}(\tilde{a}_2)\epsilon}\frac{dz}{2\pi }\frac{1-(z/z')^n}{z'-z}\\
&\times\exp\left[n\frac{{z'}^2-z^2}{2}+nix(z'-z)\right]J\omega^{(n)}(-i\tilde{a}_1z')K\omega^{(n)}(\tilde{a}_2z).
\end{split}
\end{equation}
The integrals can be interchanged as the integral over $z'$ is equivalent with a finite sum of finite terms.

Equation~\eqref{kernelprod} will be our starting point of our analysis when zooming into the origin of the macroscopic spectrum of this product matrix. This means in the present case that $\tilde{a}_1$ and $\tilde{a}_2$ are of order one.

\section{Main Result}\label{sec:result}

The aim is to take the limit $n\to\infty$ of~\eqref{kernelprod}. Due to the bounds~\eqref{lbound} and~\eqref{Kbound} of the functions $J\omega^{(n)}$ and  $K\omega^{(n)}$, respectively, the saddle points will be the same as those of the GUE, see Eq.~\eqref{saddlepoints}. The analyticity of the integrand~\eqref{kernelprod} on $\mathbb{C}\times[\mathbb{C}\setminus(i{\rm sign}(\tilde{a}_2) \mathbb{R}_+^0)]$ (Note that $z'=z$ is  not a singularity) allows us to rescale $z'\to rz'$ and $z\to rz$ with $r=|z_-|=|G(x)|$, see Eq.~\eqref{saddlepoints}. This rescaling mimics the step done at the end of subsection~\ref{sec:GUE} and it is a first step in preparing the integrand for the saddle point analysis.

\begin{figure}
\centerline{\includegraphics[width=0.675\textwidth]{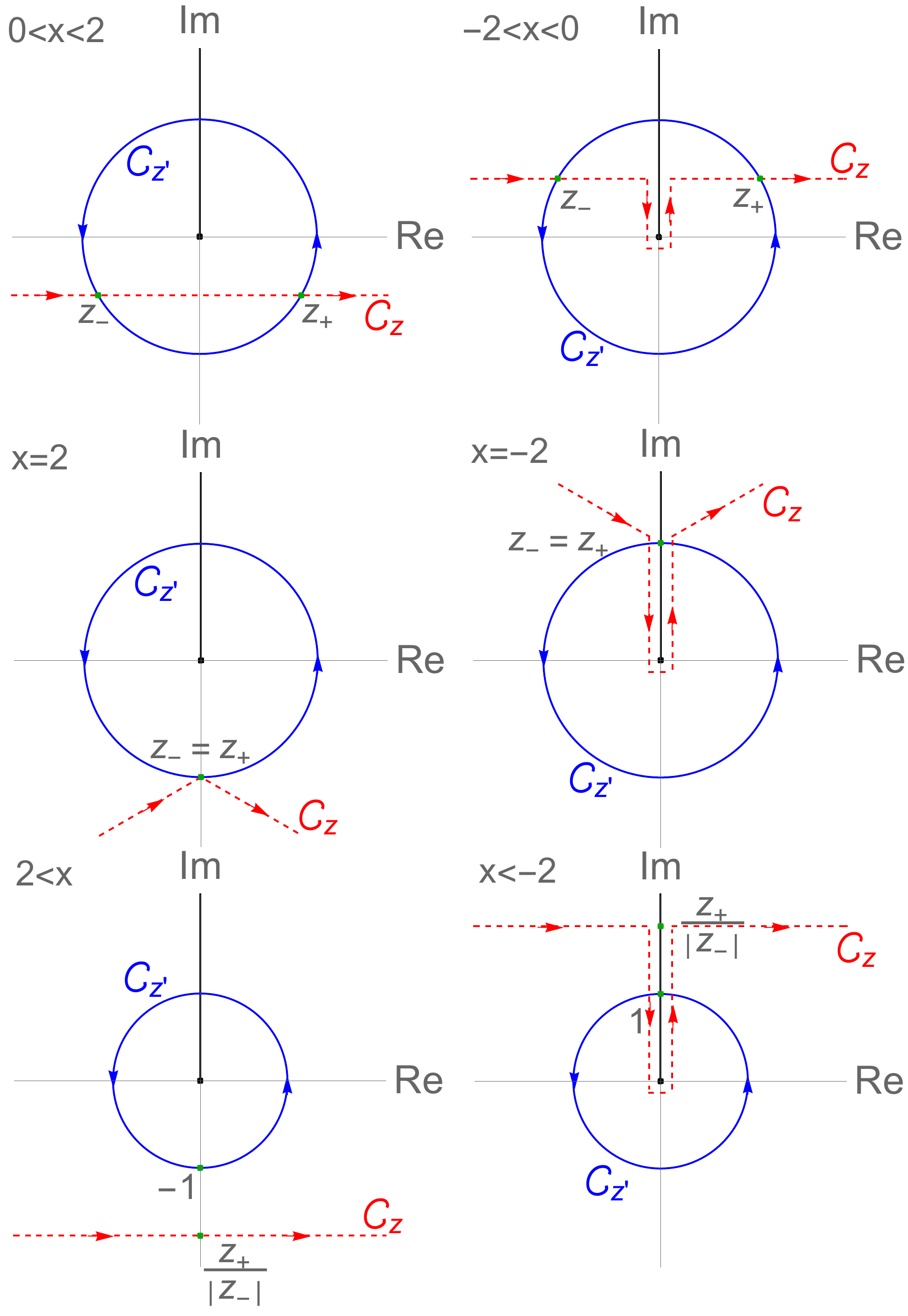}}
\caption{The two contours of the double integral~\eqref{kernelprod.b} for $\tilde{a}_2>0$. Therefore, the cut of the $K\omega^{(n)}$ is along the positive imaginary axis (black thick line with a black dot at the origin indicating the beginning of the cut). The contour $\mathcal{C}_{z'}$ (blue solid circle) is the unit circle for the $z'$-integral in all cases as well as for $x=0$. The contour $\mathcal{C}_z$ (red dashed lines) is the one of the $z$-integral where the case $x=0$ looks similar to the case $0<x<2$, with the only difference that the contour is infinitesimally close to the real axis. The saddle points~\eqref{saddlepoints} are indicated by green dots. Note, that we have rescaled the contours by $r=|z_-|$ which is smaller than $1$ for $x>2$. For the case $\tilde{a}_2<0$, the contour $\mathcal{C}_z$ as well as the cut on the imaginary axis are reflected about the real axis.}
\label{fig:countour}
\end{figure}

Let us highlight a difficulty, that we cannot easily shift the $z$-contour as we did in the GUE case, see the end of subsection~\ref{sec:GUE}, because of the cut along the half line $i{\rm sign}(\tilde{a}_2)\mathbb{R}_+^0$. Since the saddle points lie on the imaginary height $i{\rm Im}(z_+)/r\propto -i{\rm sign}(x)$, cf., Eq.~\eqref{saddlepoints}, we have to circumvent the cut whenever $\tilde{a}_2x<0$. This leads to the following contours which will be needed in Theorem~\ref{thm:main}. All the deformations are trivially guaranteed to hold due to the analyticity of the integrand in $z$ on $\mathbb{C}\setminus(i{\rm sign}(\tilde{a}_2) \mathbb{R}_+^0)$ and the Gaussian decay of the integrand at infinity whenever $\lim_{|z|\to\infty}{\rm Re}(z^2)/|z|^2>0$.

\underline{$\tilde{a}_2x>0$:}
\begin{itemize}
\item
	For $|x|\neq2$, we only shift the contour of the $z$-integral to
\begin{equation}\label{cont1}
	\mathcal{C}_z=\mathbb{R}-i\sign(x)\frac{|{\rm Im}(z_+)|}{r},
\end{equation}
see~Eq.~\eqref{saddlepoints}.
\item
	The same is essentially true for $x=0$ where it is
\begin{equation}\label{cont2}
	\mathcal{C}_z=\mathbb{R}-i\frac{{\rm sign}(\tilde{a}_2)}{n}.
\end{equation}
We have to avoid the cut still that starts at $z=0$ which amounts in the tiny shift.
\item
	For $|x|=2$ and $\tilde{a}_2x>0$, it is enough to tilt the two rays starting at $z=-i\sign(x)=-i\sign(\tilde{a}_2)$ by an angle to guarantee the integrability at the pole $z=z'$ when splitting the integral into two terms. This tilt, especially the contour, is given by
\begin{equation}\label{cont3}
	\mathcal{C}_z=\left\{-i\sign(x)+e^{-i\sign(xt)\pi/6}t|\ t\in\mathbb{R} \right\}
\end{equation}
and is allowed since one rotates the ray in a direction and through a domain where the Gaussian is still integrable.
\end{itemize}

\underline{$\tilde{a}_2x<0$:}

 In this regime, we need to split the contour into three parts, i.e.,
\begin{equation}\label{cont4}
	\mathcal{C}_z=\mathcal{C}_1\cup\mathcal{C}_2\cup\mathcal{C}_3.
\end{equation}
The part $\mathcal{C}_1$ is essentially the contours~\eqref{cont1} and~\eqref{cont3} apart from the additional condition $|{\rm Re}(z)|>1/n$ which avoids the crossing of the cut. This cut is bypassed via the contour
\begin{equation}\label{cont5}
	\mathcal{C}_2=\left\{\left.\frac{\sign(t)}{n}-i\sign(x)\left(\frac{|{\rm Im}(z_+)|}{r}t-\frac{1-|t|}{n}\right) \right|\ t\in[-1,1]\right\},
\end{equation}
which consists of two disjoint straight lines parallel to the imaginary axis. To close the contour, we need an additional section parallel to the real axis which is for all cases of $\tilde{a}_2x<0$ equal to
\begin{equation}\label{cont6}
	\mathcal{C}_3=\left[-\frac{1}{n},\frac{1}{n}\right]+i\frac{{\rm sign}(x)}{n}.
\end{equation}

All contours and all the directions of integration are shown in Fig.~\ref{fig:countour}.

One last preparation of the asymptotic analysis is the evaluation of the first term in the difference $1-(z/z')^n$ of the kernel
\begin{equation}\label{kernelprod.b}
\begin{split}
\tilde{K}_{n}(\tilde{a}_1,\tilde{a}_2)=&r\oint_{|z'|=1}\frac{dz'}{2\pi i }\int_{\mathcal{C}_z}\frac{dz}{2\pi }\frac{1-(z/z')^n}{z'-z}\exp\left[nr^2\frac{{z'}^2-z^2}{2}+nrix(z'-z)\right]\\
&\times J\omega^{(n)}(-ir\tilde{a}_1z')K\omega^{(n)}(r\tilde{a}_2z).
\end{split}
\end{equation}
Via residue theorem for the $z'$ integration, we find
\begin{equation}\label{prodkernel.c}
\tilde{K}_n(\tilde{a}_1,\tilde{a}_2)=\tilde{K}_n^{(1)}(\tilde{a}_1,\tilde{a}_2)-\tilde{K}_n^{(2)}(\tilde{a}_1,\tilde{a}_2)
\end{equation}
with the two components
\begin{equation}\label{K1}
\tilde{K}_n^{(1)}(\tilde{a}_1,\tilde{a}_2)=r\int_{\{z\in \mathcal{C}_z||z|<1\}}\frac{dz}{2\pi }J\omega^{(n)}\left(-ir\tilde{a}_1z\right)K\omega^{(n)}\left(r\tilde{a}_2z\right)
\end{equation}
and
\begin{equation}\label{K2}
\begin{split}
\tilde{K}_n^{(2)}(\tilde{a}_1,\tilde{a}_2)=&r\oint_{|z'|=1}\frac{dz'}{2\pi i }\int_{\mathcal{C}_z}\frac{dz}{2\pi }\frac{(z/z')^n}{z'-z}\exp\left[nr^2\left(\frac{{z'}^2-z^2}{2}+i\frac{x}{r}(z'-z)\right)\right]\\
&\times J\omega^{(n)}\left(-ir\tilde{a}_1z'\right)K\omega^{(n)}\left(r\tilde{a}_2z\right).
\end{split}
\end{equation}

The first part $\tilde{K}_n^{(1)}$ can be evaluated which is the following lemma about.

\begin{proposition}[Hard Edge Limit of $\tilde{K}_n^{(1)}$]\label{prop:K1}\

The limiting hard edge statistics of~\eqref{K1} is under the conditions~\ref{assumptions} given by
\begin{equation}\label{K1.a}
\begin{split}
\lim_{n\to\infty}\tilde{K}_n^{(1)}(\tilde{a}_1,\tilde{a}_2)=&\Theta(-x\tilde{a}_2)\,|{\rm Re}[ G(x)]|\,K_\infty^{(G)}( \sign(\tilde{a}_2)|{\rm Re}[ G(x)]|\tilde{a}_1,|{\rm Re}[ G(x)]\tilde{a}_2|)\\
&+\rho_{\rm GUE}(x)\int_{-1}^{1}\frac{dt}{2} J\omega^{(\infty)}\left(-i\tilde{a}_1(\pi\rho_{\rm GUE}(x)t-i{\rm Re}[ G(x)])\right)\\
&\times K\omega^{(\infty)}\left(\tilde{a}_2(\pi\rho_{\rm GUE}(x)t-i{\rm Re}[ G(x)])\right).
\end{split}
\end{equation}
for all fixed $\tilde{a}_1,\tilde{a}_2\in\mathbb{R}\setminus\{0\}$.
\end{proposition}

 The sign of $\tilde{a}_2$ can actually be absorbed in the first entry of the first term into the Green function, i.e., $\sign(\tilde{a}_2)|{\rm Re}[ G(x)]|=-{\rm Re}[ G(x)]$ because of the Heaviside step function which enforces $\sign(\tilde{a}_2)=-\sign(x)$. Indeed, this always works for any even density $\rho_{\rm GUE}(x)=\rho_{\rm GUE}(-x)$ as it is for the Wigner semicircle.

\begin{proof}
We can deform the contour in~\eqref{K1} as long as the end points (the intersections with the unit circle) remain the same and we do not cross the cut. These end points are
\begin{equation}\label{endpointstart}
z_{\rm start}=\left\{\begin{array}{cl} \displaystyle -\sqrt{1-\frac{1}{n^2}}-i\frac{\sign(\tilde{a}_2)}{n}, & x=0,\\ \displaystyle -\sqrt{1-\frac{x^2}{4}}-i\frac{x}{2}, & 0<|x|<2,\\ \displaystyle -\frac{1}{n}-i\sign(x)\sqrt{1-\frac{1}{n^2}}, & |x|\geq2\ {\rm and}\ x\tilde{a}_2<0 \end{array}\right.
\end{equation}
and $z_{\rm goal}=-\bar{z}_{\rm start}$.
We choose an $\epsilon>0$ and deform the contour for all cases such that:
\begin{enumerate}
\item 	horizontal line from $z_{\rm start}$ to $-\epsilon+i{\rm Im}(z_{\rm start})$,
\item		vertical line from $-\epsilon+i{\rm Im}(z_{\rm start})$ to $-\epsilon-i{\rm sign}(\tilde{a}_2)$,
\item 	horizontal line from $-\epsilon-i{\rm sign}(\tilde{a}_2)$ to $\epsilon-i{\rm sign}(\tilde{a}_2)$,
\item		vertical line from $\epsilon-i{\rm sign}(\tilde{a}_2)$ to $\epsilon+i{\rm Im}(z_{\rm start})$,
\item 	horizontal line from $\epsilon+i{\rm Im}(z_{\rm goal})$ to $z_{\rm goal}$.
\end{enumerate} 
When taking the limit $\epsilon\to0$, the integral over the third part vanishes since $K\omega^{(n)}$ is holomorphic in its vicinity. Therefore, we can neglect this term and concentrate on the other terms.

With the aid of Eq.~\eqref{Jomegan} and $\sign[\IM(z_{\rm start})]=-\sign(x)$ for $x>0$, the sum of the second and fourth part is essentially a difference due to the opposite direction of the integrations and they can be combined to 
\begin{equation}
\begin{split}
&\lim_{\epsilon\to0}r\int_{{\rm Im}(z_{\rm start})}^{-{\rm sign}(\tilde{a}_2)} \frac{\sign(\tilde{a}_2)idt}{2\pi} \biggl(J\omega^{(n)}\left(-ir\tilde{a}_1(it-\epsilon)\right)K\omega^{(n)}\left(r\tilde{a}_2(it-\epsilon)\right)\\
&-J\omega^{(n)}\left(-ir\tilde{a}_1(it+\epsilon)\right)K\omega^{(n)}\left(r\tilde{a}_2(it+\epsilon)\right)\biggl)\\
=&\Theta(-x\tilde{a}_2)\int_0^{|{\rm Im}(z_{\rm start})|}d(rt) J\omega^{(n)}\left(\sign(\tilde{a}_2)r\tilde{a}_1t\right)\tilde{J}\omega^{(n)}\left(r|\tilde{a}_2|t\right).
\end{split}
\end{equation}
Indeed the limit can be pushed into the integral as the integrand is dominated by an integrable function, see the bound~\eqref{Kbound} and~\eqref{Jtildebound}.
Note that for $\tilde{a}_2t<0$ the difference vanishes due to the analyticity of $K\omega^{(n)}$ in this domain.
The sign of $\tilde{a}_2$ results from the different orientations of the integrals for the two cases $\tilde{a}_2>0$ and $\tilde{a}_2<0$.
The convention of the Heaviside step function $\Theta(y)$ is chosen such that it vanishes for all arguments $y\in\mathbb{R}$ which are not positive.
Eventually we rescale $t\to |{\rm Im}(z_{\rm start})| t$ and take the limit $n\to\infty$. Then, we can be identify this part of the contour with the kernel~\eqref{hardkernPolya},
\begin{equation}\label{int1}
\begin{split}
&\lim_{n\to\infty}\lim_{\epsilon\to0}r\int_{{\rm Im}(z_{\rm start})}^{-{\rm sign}(\tilde{a}_2)} \frac{\sign(\tilde{a}_2)idt}{2\pi} \biggl(J\omega^{(n)}\left(-ir\tilde{a}_1(it-\epsilon)\right)K\omega^{(n)}\left(r\tilde{a}_2(it-\epsilon)\right)\\
&-J\omega^{(n)}\left(-ir\tilde{a}_1(it+\epsilon)\right)K\omega^{(n)}\left(r\tilde{a}_2(it+\epsilon)\right)\biggl)\\
=&\Theta(-x\tilde{a}_2)\,|{\rm Re}[ G(x)]|\,K_n^{(G)}(\sign(\tilde{a}_2) |{\rm Re}[ G(x)]|\tilde{a}_1,|{\rm Re}[ G(x)]\tilde{a}_2|).
\end{split}
\end{equation}
Here, we have exploited the relations
\begin{equation}\label{Gx.rel}
\begin{split}
{\rm Re}[ G(x)]=&-{\rm Im}( z_-)=-\lim_{n\to0}r{\rm Im}( z_{\rm start})\quad {\rm and}\\\
 \pi\rho_{\rm GUE}(x)=&|{\rm Im}[ G(x)]|=|{\rm Re}( z_-)|=\lim_{n\to0}r|{\rm Re}( z_{\rm start})|,
\end{split}
\end{equation}
where $G(x)$ has been the Green function~\eqref{GreenGUE} of the GUE and, thence, $\rho_{\rm GUE}(x)$ is the Wigner semicircle~\eqref{densGUE}.

At last we consider the sum of the first and fifth term of the contour which becomes
\begin{equation}\label{int2}
\begin{split}
&\lim_{\epsilon\to0}r\left(\int_{{\rm Re}(z_{\rm start})}^{-\epsilon}+\int_{\epsilon}^{{\rm Re}(z_{\rm goal})}\right) \frac{idt}{2\pi} J\omega^{(n)}\left(-ir\tilde{a}_1(i{\rm Im}(z_{\rm start})+t)\right)\\
&\times K\omega^{(n)}\left(r\tilde{a}_2(i{\rm Im}(z_{\rm start})+t)\right)\\
=&r\int_{-|{\rm Re}(z_{\rm start})|}^{|{\rm Re}(z_{\rm start})|}\frac{dt}{2\pi} J\omega^{(n)}\left(r\tilde{a}_1({\rm Im}(z_{\rm start})-it)\right)K\omega^{(n)}\left(r\tilde{a}_2(i{\rm Im}(z_{\rm start})+t)\right).\\
=&r|{\rm Re}(z_{\rm start})|\int_{-1}^{1}\frac{dt}{2\pi} J\omega^{(n)}\left(r\tilde{a}_1(i{\rm Im}(z_{\rm start})+|{\rm Re}(z_{\rm start})|t)\right)\\
&\times K\omega^{(n)}\left(r\tilde{a}_2({\rm Im}(z_{\rm start})-i|{\rm Re}(z_{\rm start})|t)\right).
\end{split}
\end{equation}
Note, that apart from a sign the real parts of $z_{\rm start}$ and $z_{\rm goal}$ are the same. The limit $\epsilon\to0$ can be carried out as the integrand is absolutely integrable because of the bound~\eqref{Kbound}. The same bound allows us to carry out the limit $n\to\infty$.

When combining all three contributions of $\tilde{K}_n^{(1)}(\tilde{a}_1,\tilde{a}_2)$, we arrive at our claimed statement.
\end{proof}

The following theorem, which is our main result, states that the limit~\eqref{K1.a} is actually also the limit of the full kernel $K_n$ of the product matrix $G(H-n x\eins_n)G^*$. This means the contribution of the part $K_n^{(2)}$ is always of lower order in $n$. We would like to recall that exactly the counterpart  has been tremendously important for the soft edge statistics of the GUE at $x=\pm2$.

\begin{theorem}[Hard Edge Kernel]\label{thm:main}\

The limiting hard edge statistics of the kernel~\eqref{kernelprod} for the product $G(H-n x\eins_n)G^*$ with a fixed $x\in\mathbb{R}$ and $H$ being a GUE matrix with the probability density~\eqref{probGUE} and $G$ a multiplicative P\'olya ensemble associated with the weight function $\omega^{(n)}$ that satisfies the conditions~\eqref{assump1} and~\eqref{assump2}, is given for all fixed $\tilde{a}_1,\tilde{a}_2\in\mathbb{R}\setminus\{0\}$ as follows,
\begin{equation}\label{mainresult}
\begin{split}
K_\infty(\tilde{a}_1,\tilde{a}_2)=&\lim_{n\to\infty}\tilde{K}_n(\tilde{a}_1,\tilde{a}_2)\\
=&\Theta(-{\rm Re}[ G(x)]\tilde{a}_2)\,|{\rm Re}[ G(x)]|\,K_\infty^{(G)}( -{\rm Re}[ G(x)]\tilde{a}_1,|{\rm Re}[ G(x)]\tilde{a}_2|)\\
&+\rho_{\rm GUE}(x)\int_{-1}^{1}\frac{dt}{2} J\omega^{(\infty)}\left(-i\tilde{a}_1(\pi\rho_{\rm GUE}(x)t-i{\rm Re}[ G(x)])\right)\\
&\times K\omega^{(\infty)}\left(\tilde{a}_2(\pi \rho_{\rm GUE}(x)t-i{\rm Re}[ G(x)])\right).
\end{split}
\end{equation}
The limit is pointwise.
\end{theorem}

The presentation is changed in the first few terms in Eq.~\eqref{mainresult} in comparison to~\eqref{K1.a} as we want to highlight what really enters namely the real part of the Green function. This we have seen in some exploratory computations for other ensembles than the GUE as well as when scrutinizing the proof. Equation~\eqref{K1.a} looks simpler as the sign of the real part of the Green function is directly related to $x$ via the identity~\eqref{Gx.rel} and~\eqref{endpointstart}.

Let us underline that the second term is reminiscent to the kernel derived in~\cite[Proposition 13]{FIL2018}, see also~\cite{Liu2017}, where $G$ is a product of Ginibre matrices. Thus, the theorem here is a natural, though, non-trivial generalization thereof. The kernel is different from the result~\cite[Theorem~3]{Liu2017} since in~\cite{Liu2017}, first of all, the shift of the GUE matrix $H$ is not a multiple of the identity and, secondly, most of the eigenvalues up to a finite number are symmetrically distributed at the two points $\pm1$. Therefore, the macroscopic distribution looks mostly symmetric for the model in~\cite{Liu2017}.

\begin{figure}
\centerline{\includegraphics[width=1\textwidth]{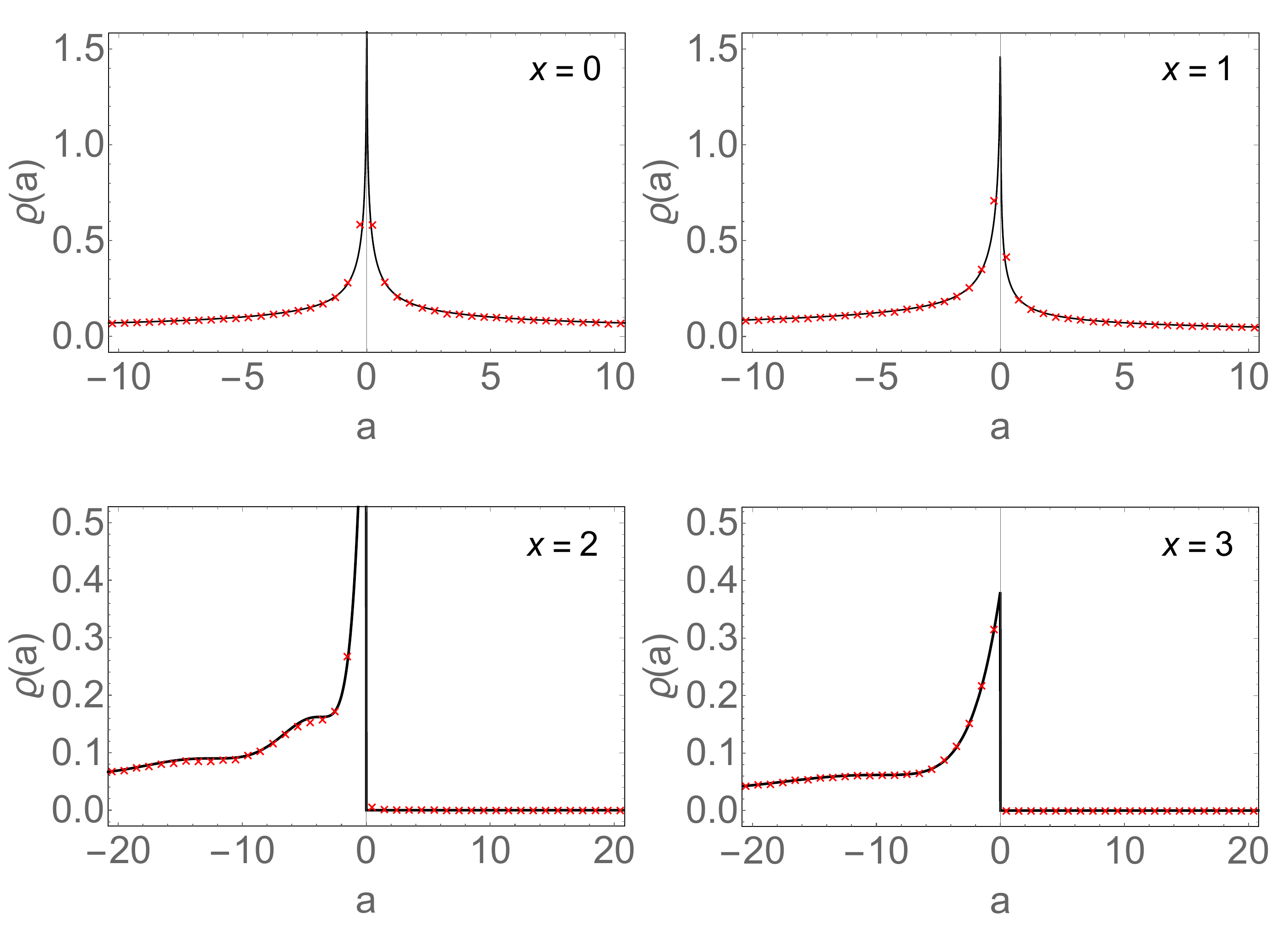}}
\caption{The microscopic level density $\rho$ at the hard edge for various shifts $x=0,1,2,3$. The analytical result~\eqref{densgin} (black solid curves) is compared with Monte Carlo simulations (red crosses). The matrix size for the numerics is chosen to be $n=100$ for $x=0,1,3$ and $n=1000$ for $x=2$. The latter is the case where $x$ sits on the soft edge of the GUE where the rate of convergence is weaker which we have also seen in the numerics due to stronger deviations and has been the reason for an increase in the matrix size. In total we have generated $10^5$ matrices for each case.}
\label{fig:hardedge}
\end{figure}

What remains to be shown is that the $K_n^{(2)}$ does indeed vanish which is proven in the next section. Here, we do not aim at an optimal rate of convergence. Nonetheless, we can say that the rate of convergence has to be at least ${\rm ln}(n)/\sqrt{n}$ and ${\rm ln}(n)/n^{1/3}$ for $x\neq\pm2$ and for $|x|=2$, respectively, cf., Sec.~\ref{sec:proof}. This is illustrated in Fig.~\ref{fig:hardedge} where we have considered the microscopic hard edge level density  of the product $G(H-n x\eins_n)G^*$ where $G$ is a Ginibre matrix drawn from the distribution
\begin{equation}
Q(G)=\pi^{-n^2}\exp[-\tr GG^*],
\end{equation}
especially we have $\omega^{(n)}(\lambda)= e^{-\lambda}$ with $\mathcal{M}\omega^{(n)}(s)=\Gamma[s]$.
This yields the three functions
\begin{equation}
\begin{split}
J\omega^{(\infty)}(z')=J_0(2\sqrt{z'}),\quad
K\omega^{(\infty)}(z)=2K_0(2\sqrt{iz}),\quad
\tilde{J}\omega^{(\infty)}(y)=J_0(2\sqrt{y}),
\end{split}
\end{equation}
where the right hand sides are the Bessel function of the first kind $J_0$ and the modified Bessel function of the  second kind $K_0$, see~\cite[Chapter~9]{ASbook} for their definitions and some of their relations.
The density is explicitly
\begin{equation}
\begin{split}\label{densgin}
\rho(a)=&K_\infty(a,a)\\
=&\Theta(-xa)|{\rm Re}[G(x)] |\int_0^1dt J_0^2\left(\sqrt{4|{\rm Re}[G(x)]a|t}\right)\\
&+\Theta(2-|x|)\frac{\sqrt{1-x^2/4}}{\pi}\int_{-1}^{1}dt I_0\left(\sqrt{2xa+i\sqrt{4-x^2}at}\right)\\
&\times K_0\left(\sqrt{2xa+i\sqrt{4-x^2}at}\right).
\end{split}
\end{equation}
The first term is the well-known Bessel kernel, only rescaled, while the second term has been derived in~\cite{FIL2018,Liu2017} for $x=0$. As already mentioned the 
analytical formula~\eqref{densgin} has been illustrated in Fig.~\ref{fig:hardedge}, where we have also included Monte Carlo simulations to underline the different order of the rates of convergence which is visible by the deviations.

Interestingly, the result~\eqref{densgin} as well as the general one~\eqref{mainresult} exhibits a nice revealing form. While the first term is the original hard edge limit of the matrix $GG^*$, the second term is essentially (apart from subtle details) the transformation formula~\eqref{kerneltrafo} applied to  the sine-kernel
\begin{equation}
\frac{\sin[\tilde{a}_1/\tilde{z}-\tilde{a}_2/a]}{\tilde{a}_1/\tilde{z}-\tilde{a}_2/a}=\frac{1}{2}\int_{-1}^1dt\exp\left[i\left(\frac{\tilde{a}_1}{\tilde{z}}-\frac{\tilde{a}_2}{a}\right)t\right],
\end{equation}
cf., the definitions~\eqref{Iomegan} and~\eqref{Komegan}.
Indeed when replacing the two Bessel functions in the second integral by the Fourier exponentials $\exp[2x\tilde{a}_1+i\sqrt{4-x^2}\tilde{a}_1t]$ and $\exp[-2x\tilde{a}_2-i\sqrt{4-x^2}\tilde{a}_2t]$, the kernel reduces to the original sine-kernel result.  Thus, these two terms seem to be a very natural decomposition.

\section{Proof of Theorem~\ref{thm:main}}\label{sec:proof}

The aim of this section is to show that the two-fold integral
\begin{equation}\label{Iint}
\begin{split}
\mathcal{I}=&r\oint_{|z'|=1}\frac{dz'}{2\pi i }\int_{\mathcal{C}_z}\frac{dz}{2\pi }\frac{(z/z')^n}{z'-z}\exp\left[nr^2\left(\frac{{z'}^2-z^2}{2}+i\frac{x}{r}(z'-z)\right)\right]\\
&\times J\omega^{(n)}\left(-i\tilde{a}_1rz'\right)K\omega^{(n)}\left(\tilde{a}_2rz\right)
\end{split}
\end{equation}
vanishes for $n\to\infty$ regardless which case of $x$ and $a_2\neq0$ we consider.  Let us emphasize that we do not aim at the optimal rate of convergence of this integral since that one strongly depends on the considered weight $\omega^{(n)}$.

The two poles at $z=z'=z_\pm$ for $|x|<2$  and at $z=z'=-i\sign(x)e^{\pm i\epsilon}$ for $|x|\geq2$ are integrable as it is jointly a two-dimensional integral. When $z$ approaches the circle from the inside it always crosses the $z'$-contour at a non-vanishing angle $0<\theta\leq \pi/2$. Thence,  there is an open neighbourhood at $z=z'$ which can be mapped to the complex half plane where the measure $dz' dz/|z'-z|$ becomes $d^2\tilde{z}/|\tilde{z}|= dR d\phi$ with $\tilde{z}=Re^{i\phi}$. Similarly this is the case when $z$ comes from the outside, though then we have also the case that no crossing is going to happen for $\tilde{a}_2x>0$ and $|x|>2$. In particular we make use of the estimate
\begin{equation}\label{estimate.ratio}
\frac{|1-r+i(\varphi-\vartheta)/\pi|}{|e^{i\varphi}-re^{i\vartheta}|}=\sqrt{\frac{([\varphi-\vartheta]/\pi)^2+(1-r)^2}{\sin^2[\varphi-\vartheta]+(\cos[\varphi-\vartheta]-r)^2}}\leq \sqrt{2}
\end{equation}
for all $r\in\mathbb{R}_+$ and $(\varphi-\vartheta)\in]-\pi,\pi[$. This result can be obtained by noticing that the ratio has two extrema in $r$ (first derivative reduces to a quadratic equation) while one needs to be careful that $r$ has to be positive. Only one of the two extrema of $r$ is positive which is even a minimum. Therefore, one only needs to consider the extrema $r\to0$ and $r\to\infty$ where the first is evidently the largest of the two. Moreover, the maximum in $r$ achieves its the supremum in the limit $|\varphi-\vartheta|\to\pi$.

Additionally, the functions $ J\omega^{(n)}\left(-ia_1rz'\right)$ and $K\omega^{(n)}\left(a_2rz\right)$ are bounded by exponential functions, see~\eqref{lbound} and~\eqref{Kbound}, which are sub-leading terms compared to the Gaussian bound at infinity in the integrand. Therefore, we are left by estimating the terms
\begin{equation}
\frac{F(z)}{F(z')}=\exp\left[nr^2\left(\frac{{z'}^2}{2}+i\frac{x}{r}z'\right)-nr^2\left(\frac{{z}^2}{2}+i\frac{x}{r}z\right)\right]\left(\frac{z}{z'}\right)^n.
\end{equation}

First, we find an upper bound for the term corresponding to $z'=e^{\imath\varphi}$ which is
\begin{equation}\label{zprime.est}
\begin{split}
\left|\frac{1}{F(z')}\right|=&\exp\left[nr^2\left(\frac{\cos(2\varphi)}{2}-\frac{x}{r}\sin\varphi\right)\right]\leq\left\{\begin{array}{cl} \displaystyle \exp\left[n\frac{2+x^2}{4}\right], & |x|\leq 2, \\ \displaystyle \exp\left[n\frac{2+|z_-|^2}{2}\right], & |x|\geq2.\end{array}\right.
\end{split}
\end{equation}
The maximum is reached at $z'=z_\pm$ for $|x|<2$ and at $z'=-i\sign(x)$ for $|x|\geq2$. That has been the simplest part of the proof.  Now we turn to $F(z)$.

Before we go over to the case discussion, let us mention a peculiarity of the notation that has been introduced to focus on the $n$-dependence and not the detailed expression of the proportionality constants.
When we write for two expressions $\mathcal{E}_1$ and $\mathcal{E}_2$
\begin{equation}
\mathcal{E}_1\leq{\rm const.}\ \mathcal{E}_2\quad{\rm or}\quad \mathcal{E}_1={\rm const.}\ \mathcal{E}_2
\end{equation}
it means that there is an $n$-independent constant, represented by ${\rm const.}$, so that these relations hold. Certainly, these constants can dependent on $\tilde{a}_1$, $\tilde{a}_2$ and/or $x$. We do not exclude that. However, we want to show that the double contour integral~\eqref{Iint} vanishes for $n\to\infty$. That is the primary goal.

\subsection{The Case $x=0$}

For $x=0$, we have the parametrisation $z=t-i\,\sign(\tilde{a}_2)/n$ which yields the estimate
\begin{equation}\label{Fx00}
\begin{split}
\left|F\left(t-i\frac{\sign(\tilde{a}_2)}{n}\right)\right|=&\exp\left[-n\left(\frac{t^2-1/n^2}{2}\right)\right]\left(t^2+\frac{1}{n^2}\right)^{n/2}\leq \exp\left[-\frac{n^2-2}{2n}\right],
\end{split}
\end{equation}
because the maximum is reached at $t^2=1-1/n^2$ while $t=0$ is a local minimum. We combine this with the fact that $ J\omega^{(n)}\left(-i\tilde{a}_1e^{i\varphi}\right)$ is bounded by  $\widetilde{C}e^{|\tilde{a}_1|}$ and $K\omega^{(n)}\left(\tilde{a}_2t-i|\tilde{a}_2|/n\right)$ by $c/(|\tilde{a}_2|\sqrt{t^2+1/n^2})$ for some constant $c>0$. Hence, we have in total
\begin{equation}
\begin{split}
|\mathcal{I}|\leq& {\rm const.}\int_{-\pi}^\pi\frac{d\varphi}{2\pi }\int_{-\infty}^\infty\frac{dt}{2\pi }\frac{\left(t^2+1/n^2\right)^{(n-1)/2}}{|e^{i\varphi}-t+i\,\sign(\tilde{a}_2)/n|}\exp\left[-\frac{n^2(t^2-1)-1}{2n}\right]\\
\leq& {\rm const.}\int_{-\pi}^\pi\frac{d\varphi}{2\pi }\int_{0}^\infty\frac{dt}{2\pi }\sum_{L=\pm1}\frac{\left(t^2+1/n^2\right)^{(n-1)/2}}{\sqrt{({\rm mod}_{]-\pi,\pi]}[\varphi-\vartheta(Lt)])^2/\pi^2+(1-\sqrt{t^2+1/n^2})^2}}\\
&\times\exp\left[-\frac{n^2(t^2-1)-1}{2n}\right].
\end{split}
\end{equation}
The second inequality is a consequence of Eq.~\eqref{estimate.ratio}, where we used the abbreviation $Lt-i\,\sign(\tilde{a}_2)/n=\sqrt{t^2+1/n^2}e^{i\vartheta(Lt)}$ and a shifted modulo ${\rm mod}_{]-\pi,\pi]}[\phi]=\tilde{\phi}\in]-\pi,\pi]$ for any $\phi=\tilde{\phi}+2\pi l$ with $l\in\mathbb{Z}$. The angle $\vartheta(Lt)$ can be absorbed by $\varphi\to\varphi +\vartheta(Lt)$ for each of the two summands, separately, leading to a factor of $2$ and the integral over $\varphi$ can be carried out exactly,
\begin{equation}
\begin{split}
|\mathcal{I}|\leq&{\rm const.}\int_{0}^\infty\frac{dt}{2\pi }{\rm arcsinh}\left[|1-\sqrt{t^2+1/n^2}|^{-1}\right]\exp\left[-\frac{n^2(t^2-1)-1}{2n}\right]\\
&\times\left(t^2+\frac{1}{n^2}\right)^{(n-1)/2}.
\end{split}
\end{equation}

In the next step, we substitute $t=\sqrt{(1+\delta t/\sqrt{n})^2-1/n^2}$ with $\delta t\in[-\sqrt{n}+1/\sqrt{n},\infty[$, and we obtain
\begin{equation}
\begin{split}
|\mathcal{I}|\leq& \frac{{\rm const.}}{\sqrt{n}}\int_{-\sqrt{n}+1/\sqrt{n}}^\infty \hspace*{-0.5cm}d\delta t\frac{{\rm arcsinh}\left[\sqrt{n}/|\delta t|\right]\left(1+\delta t/\sqrt{n}\right)^{n}}{2\pi \sqrt{(1+\delta t/\sqrt{n})^2-1/n^2}}\exp\left[-\frac{\delta t^2}{2}-\sqrt{n}\delta t+\frac{1}{n}\right]\\
\leq&\frac{{\rm const.}}{\sqrt{n}}\int_{-\sqrt{n}+1/\sqrt{n}}^\infty\hspace*{-0.5cm}d\delta t\frac{{\rm arcsinh}(1)+\left|{\rm ln}\left[\sqrt{n}/|\delta t|\right]\right|}{2\pi \sqrt{(1+\delta t/\sqrt{n})^2-1/n^2}}\exp\left[-\frac{\delta t^2}{2}\right]
\end{split}
\end{equation}
for all $n\geq1$. Here, we employed the estimates ${\rm arcsinh}(|x|)\leq {\rm arcsinh}(1)+|{\rm ln}(|x|)|$ and $(1+x)^ne^{-n x}\leq 1$ for all $x\geq-1$. The integral is then split into the intervals $]-\sqrt{n}+1/\sqrt{n},-\sqrt{n}/2]$ and $]-\sqrt{n}/2,\infty[$ so that we can either estimate the Gaussian or the denominator. Thus, we arrive at
\begin{equation}
\begin{split}
|\mathcal{I}|\leq&  {\rm const.}\left[\int_{1/n}^{1/2}\frac{d\delta t}{2\pi \sqrt{\delta t^2-1/n^2}}\left({\rm arcsinh}(1)-{\rm ln}\left[|\delta t-1|\right]\right)\exp\left[-\frac{n}{8}\right]\right.\\
&\left.+\frac{1}{\sqrt{n}}\int_{-\sqrt{n}/2}^\infty\frac{d\delta t}{2\pi \sqrt{1/4-1/n^2}}\left({\rm arcsinh}(1)+{\rm ln}\left[\frac{\sqrt{n}}{|\delta t|}\right]\right)\exp\left[-\frac{\delta t^2}{2}\right]\right]
\end{split}
\end{equation}
which is now valid only for $n> 2$. Notice, that we made a change, $\delta t\to \sqrt{n}(\delta t-1)$, in the first integral. Both integrals are finite but the second is logarithmically growing in $n$. Hence, the first is exponentially suppressed by the prefactor $\exp[-n/8]$ and the second one behaves like ${\rm ln}(n)/\sqrt{n}$ for large $n$. Therefore, we arrive at the final estimate
\begin{equation}
|\mathcal{I}|\leq{\rm const.}\ \frac{{\rm ln}(n)}{\sqrt{n}}
\end{equation}
with some numerical constant. It evidently vanishes for each fixed $\tilde{a}_1\in\mathbb{C}$ and $\tilde{a}_2\in\mathbb{C}\setminus\{0\}$ when taking $n\to\infty$.

\subsection{The Case $0<|x|<2$}

For $|x|< 2$, the parametrization $z=t-ix/2$ with real $t$ leads to
\begin{equation}\label{Fx}
\begin{split}
\left|F\left(t-i\frac{x}{2}\right)\right|=&\exp\left[-n\left(\frac{t^2-x^2/4}{2}+\frac{x^2}{2}\right)\right]\left(t^2+\frac{x^2}{4}\right)^{n/2}\leq \exp\left[-n\frac{2+x^2}{4}\right],
\end{split}
\end{equation}
because the maximum is acquired at $t^2=1-x^2/4$ while $t=0$ is again a minimum. Now we have to distinguish the two cases of $\tilde{a}_2x>0$ or $\tilde{a}_2x<0$.

In the first case, we do not need to make a detour and can choose $t\in\mathbb{R}$. Then everything works along the same lines as for the case $x=0$ with the difference that the bound of $K\omega^{(n)}\left(\tilde{a}_2(t-ix/2)\right)$ simplifies to $2c/(\tilde{a}_2 x)$, since we are still on the  half plane ${\rm Im}(\tilde{a}_2z)<0$. Thence, we have
\begin{equation}
\begin{split}
|\mathcal{I}|\leq& {\rm const.}\int_{-\pi}^\pi\frac{d\varphi}{2\pi }\int_{-\infty}^\infty\frac{dt}{2\pi }\frac{1}{|e^{i\varphi}-t+i\,x/2|}\exp\left[-n\left(\frac{t^2-1}{2}+\frac{x^2}{8}\right)\right]\left(t^2+\frac{x^2}{4}\right)^{n/2}\\
\leq&{\rm const.}\int_{0}^\infty\frac{dt}{2\pi }{\rm arcsinh}\left(\left|1-\sqrt{x^2/4+t^2}\right|^{-1}\right)\exp\left[-n\left(\frac{t^2-1}{2}+\frac{x^2}{8}\right)\right]\\
&\times\left(t^2+\frac{x^2}{4}\right)^{n/2}.
\end{split}
\end{equation}
In the last line, we again applied Eq.~\eqref{estimate.ratio} and carried out the integration over $\varphi$. As before we substitute $t=\sqrt{(1+\delta t/\sqrt{n})^2-x^2/4}$ with $\delta t\in](-1+|x|/2)\sqrt{n},\infty[$ and obtain
\begin{equation}\label{h.eq1}
\begin{split}
|\mathcal{I}|\leq&\frac{{\rm const.}}{\sqrt{n}}\int_{0}^\infty\frac{d\delta t}{2\pi \sqrt{(1+\delta t/\sqrt{n})^2-x^2/4}}{\rm arcsinh}\left(\frac{\sqrt{n}}{|\delta t|}\right)\exp\left[-\frac{\delta t^2}{2} -\sqrt{n}\delta t\right]\\
&\times\left(1+\frac{\delta t}{\sqrt{n}}\right)^{n+1}\\
\leq &\frac{{\rm const.}}{\sqrt{n}}\int_{(-1+|x|/2)\sqrt{n}}^\infty\frac{d\delta t}{2\pi \sqrt{(1+\delta t/\sqrt{n})^2-x^2/4}}\left({\rm arcsinh}(1)+\left|{\rm ln}\left[\frac{\sqrt{n}}{|\delta t|}\right]\right|\right)\\
&\times\exp\left[-\frac{\delta t^2}{2} \right]\left(1+\frac{\delta t}{\sqrt{n}}\right).
\end{split}
\end{equation}
This integral is not that different from the one in the case $x=0$ case, so that after splitting the integral into two parts over the intervals $](-1+|x|/2)\sqrt{n},(-1+|x|/2)\sqrt{n}/2[$  and $](-1+|x|/2)\sqrt{n}/2,\infty[$ the bound is again given by
\begin{equation}\label{boundin}
\begin{split}
|\mathcal{I}|\leq&{\rm const.}\frac{{\rm ln}(n)}{\sqrt{n}}
\end{split}
\end{equation}
for all $n>2$, which obviously vanishes for $n\to\infty$.

For $\tilde{a}_2x<0$, we only consider large enough $n\in\mathbb{N}$ such that ${\rm ln}(2/n^2)+1+x^2/2<0$ which is  indeed given for $n>\sqrt{2}e^{1/2+x^2/4}$. Then, we have ${\rm ln}(t^2+1/n^2)+1+x^2/2<0$ for all $|t|<1/n$. We will need this later on.

The contour is decomposed into three parts $\mathcal{C}_z=\mathcal{C}_1\cup\mathcal{C}_2\cup\mathcal{C}_3$.
In the first part $\mathcal{C}_1$, we integrate over $z=t-ix/2$ with $|t|>1/n$ and $t\in\mathbb{R}$. This contribution is coined by $\mathcal{I}_{\mathcal{C}_1}$; the subscript reflects its association of the contour. We consider almost the same integral as for $\tilde{a}_2x>0$ apart from the estimate $|K\omega^{(n)}(r\tilde{a}_2z)|\leq{\rm const.}(|z|+1)/|z|\, \exp[\alpha|z|]$, where $\alpha>0$ is $n$-independent, which follows from~\eqref{Kbound}. In the variable $\delta t$, the  factor $(|z|+1)\exp[\alpha|z|]$ reads $(2+\delta t/\sqrt{n})\exp[\alpha(1+\delta t/\sqrt{n})]$ which is dominated by the Gaussian in Eq.~\eqref{h.eq1}. Thence, nothing changes with the estimation and $\mathcal{I}_{\mathcal{C}_1}$ is still bounded by ${\rm ln}(n)/\sqrt{n}$ up to a constant.

 The other parallel part $\mathcal{C}_3$ with $z=t-i\,\sign(\tilde{a}_2)/n$ for $|t|<1/n$ is given by
\begin{equation}
\begin{split}
|\mathcal{I}_{\mathcal{C}_3}|\leq& {\rm const.}\int_{-\pi}^\pi\frac{d\varphi}{2\pi }\int_{-1/n}^{1/n}\frac{dt}{2\pi }\frac{\left(t^2+1/n^2\right)^{(n-1)/2}}{|e^{i\varphi}-t+i\,\sign(\tilde{a}_2)/n|}\\
&\times\exp\left[-\frac{n^2(t^2-1)-1}{2n}+n\frac{x^2}{4}+|x|\right],
\end{split}
\end{equation}
where we have exploited Eq.~\eqref{Fx00} and needed again the bound of $K\omega^{(n)}\left(\tilde{a}_2(t-i\frac{\sign(\tilde{a}_2)}{n})\right)$ from the case $x=0$. The denominator $|e^{i\varphi}-t+i\,\sign(\tilde{a}_2)/n|$ stays this time away from zero by the constant $(1-\sqrt{2}/n)/\sqrt{2}\geq 1/2>0$ for all $n\geq3$, leading to the simplification
\begin{equation}\label{apprxl2.a}
\begin{split}
|\mathcal{I}_{\mathcal{C}_3}|\leq& {\rm const.}\int_{0}^{1/n}\frac{dt}{2\pi }\left(t^2+\frac{1}{n^2}\right)^{(n-1)/2}\exp\left[-\frac{n^2(t^2-1)-1}{2n}+n\frac{x^2}{4}+|x|\right]\\
\leq& {\rm const.}\int_{0}^{1/n}\frac{dt}{2\pi }\left(t^2+\frac{1}{n^2}\right)^{-1/2}\left(\frac{2}{n^2}\right)^{n/2}\exp\left[\frac{n}{2}+n\frac{x^2}{4}+|x|\right]\\
\leq& {\rm const.}\left(\frac{2}{n^2}\right)^{n/2}\exp\left[\frac{n}{2}+n\frac{x^2}{4}\right].
\end{split}
\end{equation}
Exactly here, we have exploited the fact $-t^2+{\rm ln}(t^2+1/n^2)+1+x^2/2\leq -1/n^2+{\rm ln}(2/n^2)+1+x^2/2<0$ because the maxima in $t$, namely $t^2=1-1/n^2$ lie outside of the interval. This also tells us that this term vanishes exponentially, like $e^{-n\ {\rm ln}(n)+\gamma n}$ for some $\gamma\in\mathbb{R}$, in the limit $n\to\infty$.

At last, we consider the integration parallel to the imaginary axis, i.e., $\mathcal{C}_2=\{z=\sign(t)/n-i\sign(x) [|tx|/2-(1-|t|)/n] | t\in[-1,1]\}$, where we have
\begin{equation}\label{Fxi}
\begin{split}
\left|F\left(\frac{\sign(t)}{n}-i\sign(x) \lambda\right)\right|=&\exp\left[-n\left(\frac{1/n^2-\lambda^2}{2}+|x|\lambda\right)\right]\left(\frac{1}{n^2}+\lambda^2\right)^{n/2}
\end{split}
\end{equation}
with $\lambda=|tx|/2-(1-|t|)/n\in[-1/n,|x|/2]$. The denominator $1/|e^{i\varphi}-z|$ is again bounded from above by a constant $\sqrt{2}/|1-\sqrt{1/n^2+x^2/4}|\leq\sqrt{2}/(1-\sqrt{1/n_0^2+x^2/4})<\infty$ for all $n>n_0>1/\sqrt{1-x^2/4}$. The maximum $K\omega^{(n)}\left(\tilde{a}_2(\frac{\sign(t)}{n}-i\lambda)\right)$ is given either by the point closest to the origin or at its extremal points. When $n$ is large enough, this maximum is given by $\lambda=0$ so that  we can take the bound $cn e^{\alpha|\tilde{a}_2|/n}/|\tilde{a}_2|\leq cn e^{\alpha|\tilde{a}_2|}/|\tilde{a}_2|$. Then, we obtain
\begin{equation}
\begin{split}
|\mathcal{I}_{\mathcal{C}_2}|\leq& {\rm const.}\ n\int_{-1/n}^{|x|/2}\frac{d\lambda\ \left(\frac{1}{n^2}+\lambda^2\right)^{n/2}}{2\pi (|x|/2+1/n)|\tilde{a}_2|}\exp\left[-n\left(\frac{1/n^2-\lambda^2}{2}+|x|\lambda-\frac{2+x^2}{4}\right)\right],
\end{split}
\end{equation}
where each of the two components of $\mathcal{C}_2$ yield the same. The integral over the interval $[-1/n,0]$ is again exponentially bounded by $\exp[n/2({\rm ln}(2/n^2)+1+x^2/2)]$ for all $n\geq \sqrt{2}e^{1/2+x^2/4}$. The integral over the positive part of the interval can be estimated from above by the upper bound at $\lambda=|x|/2$ since the integrand is an increasing function for $0<\lambda<|x|/2$, when $n$ is big enough, as can be readily checked by taking the derivative in $\lambda$ of the exponent. The leading term in the exponent is then  $n/2(1-x^2/4-1/n^2+{\rm ln}[x^2/4+1/n^2])$ which is indeed always negative apart from the point $x^2/4=1-1/n^2$. Therefore, also this part is exponentially bounded when $n$ is large enough and $|x|$ stays away from $2$.

Summarizing, also for $\tilde{a}_2x<0$, we obtain the bound~\eqref{boundin} for the contour integral $\mathcal{I}$.

\subsection{The Case $|x|=2$}

Also for this case, we start with the simpler version of $\tilde{a}_2x>0$. The parametrization of the $z$-integral is given by $z=-i\sign(x)+e^{-i\sign(x t)\pi/6}t$ with $t\in\mathbb{R}$. This leads to
\begin{equation}\label{Fxon}
\begin{split}
\left|F\left(-i\sign(x)+e^{-i\sign(x t)\pi/6}t\right)\right|=&\exp\left[-n\frac{t^2+2|t|+6}{4}\right]\left(t^2+|t|+1\right)^{n/2}.
\end{split}
\end{equation}
The term $1/|e^{i\varphi}-z|$ has the upper bound $\sqrt{2}/\sqrt{(\varphi-\vartheta(t))^2/\pi^2+(\sqrt{1+|t|+t^2}-1)^2}$ for $[\varphi-\vartheta(t)]\in]-\pi,\pi[$ and $K\omega^{(n)}(\tilde{a}_2(-i\sign(x)+e^{-i\sign(x t)\pi/6}t))$ can be replaced by $c/(|\tilde{a}_2|\sqrt{1+|t|+t^2})\leq c/|\tilde{a}_2|$. Note that both rays, $t>0$ and $t<0$, can be estimated by the same bound after splitting the integral in its two components and then integrating over $\varphi$. Thus, we find
\begin{equation}
\begin{split}
|\mathcal{I}|\leq&{\rm const.}\int_{0}^\infty\frac{dt}{2\pi }{\rm arcsinh}\left(\frac{1}{\sqrt{1+t+t^2}-1}\right)\exp\left[-n\frac{t^2+2t}{4}\right]\left(t^2+t+1\right)^{n/2}.
\end{split}
\end{equation}
In the next step, we substitute $t=\sqrt{e^{\delta t/n^{1/3}}-3/4}-1/2$ with $\delta t>0$ and obtain
\begin{equation}
\begin{split}
|\mathcal{I}|\leq&\frac{{\rm const.}}{n^{1/3}}\int_{0}^\infty\frac{d\delta t}{2\pi \sqrt{4e^{\delta t/n^{1/3}}-3}}{\rm arcsinh}\left(\frac{1}{e^{\delta t/(2n^{1/3})}-1}\right)\\
&\times\exp\left[-\frac{n}{4}\left(e^{\delta t/n^{1/3}}+\sqrt{e^{\delta t/n^{1/3}}-\frac{3}{4}}-\frac{3}{2}\right)+\left(1+\frac{n}{2}\right)\frac{\delta t}{n^{1/3}}\right].
\end{split}
\end{equation}
The denominator $1/ \sqrt{4e^{\delta t/n^{1/3}}-3}$ can be bounded by $1$ while the ${\rm arcsinh}$ function has the upper bound
\begin{equation}
\begin{split}
{\rm arcsinh}\left(\frac{1}{e^{\delta t/(2n^{1/3})}-1}\right)\leq& {\rm arcsinh}\left(1\right)+\left|{\rm ln}\left(e^{\delta t/(2n^{1/3})}-1\right)\right|\\
\leq&{\rm arcsinh}\left(1\right)+\left|{\rm ln}\left(\frac{\delta t}{2n^{1/3}}\right)\right|+\frac{\delta t}{2n^{1/3}}.
\end{split}
\end{equation}
The exponent can be approximated from above as
\begin{equation}
\begin{split}
&-\frac{n}{4}\left(e^{\delta t/n^{1/3}}+\sqrt{e^{\delta t/n^{1/3}}-\frac{3}{4}}-\frac{3}{2}\right)+\left(1+\frac{n}{2}\right)\frac{\delta t}{n^{1/3}}\\
\leq&-\frac{n}{4}\left(1+\frac{\delta t}{n^{1/3}}+\frac{\delta t^2}{2n^{2/3}}+\frac{\delta t^3}{6n}+\sqrt{1+\frac{\delta t}{n^{1/3}}+\frac{\delta t^2}{2n^{2/3}}-\frac{3}{4}}-\frac{3}{2}\right)+\left(1+\frac{n}{2}\right)\frac{\delta t}{n^{1/3}}\\
\leq&-\frac{\delta t^3}{24}+\frac{\delta t}{n^{1/3}}-\frac{n}{4}\left(\frac{\delta t^2}{2n^{2/3}}-\frac{\delta t}{n^{1/3}}+\sqrt{\frac{\delta t^2}{2n^{2/3}}+\frac{\delta t}{n^{1/3}}+\frac{1}{4}}-\frac{1}{2}\right)\\
\leq&-\frac{\delta t^3}{24}+\frac{\delta t}{n^{1/3}}.
\end{split}
\end{equation}
The approximation of the exponential functions by a finite sum is allowed since its argument is positive ($\delta t>0$) and, hence, the sum is smaller than the full series. In the last step, we used the fact that the expression in the parenthesis is positive for any $\delta t>0$ which can be found by differentiating in $\tilde{t}=\delta t/n^{1/3}+1$. Its second derivative is positive for any $\tilde{t}>1$ so that the first derivative is strictly increasing and, therefore, positive, too. Thus, the minimum is taken at $\tilde{t}=1$ or, equivalently, $\delta t=0$ which vanishes.

Collecting the discussion above, we eventually have
\begin{equation}
\begin{split}
|\mathcal{I}|\leq&\frac{{\rm const.}}{n^{1/3}}\int_{0}^\infty\frac{d\delta t}{2\pi }\left({\rm arcsinh}\left(1\right)+\left|{\rm ln}\left(\frac{\delta t}{2n^{1/3}}\right)\right|+\frac{\delta t}{2n^{1/3}}\right)\exp\left[-\frac{\delta t^3}{24}+\frac{\delta t}{n^{1/3}}\right]\\
\leq&{\rm const.}\frac{{\rm ln}(n)}{n^{1/3}}.
\end{split}
\end{equation}
This bound vanishes for each fixed $\tilde{a}_1$ and $\tilde{a}_2$.

Let us turn to the case $\tilde{a}_2x<0$. The contribution of the contour $\mathcal{C}_1=\{z=-i\sign(x)+e^{-i\sign(x t)\pi/6}t|\sqrt{3}|t|/2>1/n\}$ can be computed in the same way as for the case $\tilde{a}_2x>0$, despite the fact that the function $K\omega^{(n)}$ yields an additional exponential term $c (|\tilde{a}_2|\sqrt{1+|t|+t^2}+1)/(|\tilde{a}_2|\sqrt{1+|t|+t^2})\,\exp[\alpha(|\tilde{a}_2|\sqrt{1+|t|+t^2})]$ $\leq c(1+1/|\tilde{a}_2|)\exp[\alpha(|\tilde{a}_2|\sqrt{1+|t|+t^2})]$, cf. Eq.~\eqref{Kbound}. To see this one can show for $t^2+|t|+1=e^{\delta t/n^{1/3}}>1/n^2+1/n+1$ that the term $-n/4 e^{\delta t/n^{1/3}}+\alpha|\tilde{a}_2| e^{\delta t/(2n^{1/3})}$ has a Taylor series with only negative coefficients for all  $n>4\alpha|\tilde{a}_2|$. Therefore, we can employ
\begin{equation}
\begin{split}
-\frac{n}{4}e^{\delta t/n^{1/3}}+\alpha|\tilde{a}_2| e^{\delta t/(2n^{1/3})}\leq& -\frac{n}{4}\left(1+\frac{\delta t}{n^{1/3}}+\frac{\delta t^2}{2n^{2/3}}+\frac{\delta t^3}{6n}\right)\\
&+\alpha|\tilde{a}_2| \left(1+\frac{\delta t}{2n^{1/3}}+\frac{\delta t^2}{8n^{2/3}}+\frac{\delta t^3}{48n}\right),
\end{split}
\end{equation}
which leads us onto the old track where we found the upper bound of ${\rm ln}(n)/n^{1/3}$.

The estimate for $\mathcal{C}_3=\{z=t-i\,\sign(\tilde{a}_2)/n|\ t\in[-1/n,1/n]\}$ does not differ at all from the case $|x|<2$. Therefore, we have with Eq.~\eqref{apprxl2.a}
\begin{equation}\label{apprx2.a}
\begin{split}
|\mathcal{I}_{\mathcal{C}_3}|\leq&  {\rm const.}\left(\frac{2}{n^2}\right)^{n/2}\exp\left[\frac{n}{2}+n+2\right]
\end{split}
\end{equation}
for $-t^2+{\rm ln}(t^2+1/n^2)+3\leq -1/n^2+{\rm ln}(2/n^2)+3<0$, meaning $n$ has to be big enough. The bound~\eqref{apprx2.a} vanishes due its leading behaviour $e^{-n\,{\rm ln}(n)}$.

Finally, we consider the part $\mathcal{C}_2=\{z=\sign(t)/n-i\sign(x) [(1+1/(\sqrt{3}n))|t|-(1-|t|)/n]|\ t\in[-1,1]\}$.
Note that Eq.~\eqref{Fxi} is still valid only that we have now $\lambda=(1+1/(\sqrt{3}n))|t|-(1-|t|)/n$ with $\lambda\in[-1/n,1+1/(\sqrt{3}n)]$. Since we cross with the $z$-integral the circle described by $z'=e^{i\varphi}$,  we cannot easily replace $1/|e^{i\varphi}-z|$ by a constant as we could for the case $|x|<2$. We need to employ Eq.~\eqref{estimate.ratio} and integrate for each of the two disjoint parts of $\mathcal{C}_2$ over $\varphi$. Exchanging $K\omega^{(n)}$ by the bound $c/(|\tilde{a}_2|\sqrt{\lambda^2+1/n^2}) $, where we absorb the constant bound for the exponential term $(|z|+1)\exp[\alpha|z|]\leq (\sqrt{5}|\tilde{a}_2|+1)\exp[\sqrt{5}\alpha |\tilde{a}_2|]$ in $c$. This leads us to
\begin{equation}\label{apprx2.b}
\begin{split}
|\mathcal{I}_{\mathcal{C}_2}|\leq& {\rm const.}\biggl(\int_{-1/n}^{1/2}+\int_{1/2}^{1}+\int_{1}^{1+1/(\sqrt{3}n)}\biggl)\frac{d\lambda\ \left(\frac{1}{n^2}+\lambda^2\right)^{(n-1)/2}}{2\pi(1+1/(\sqrt{3}n)+1/n)}\\
&\times{\rm arcsinh}\left(\frac{1}{|1-\sqrt{\lambda^2+1/n^2}|}\right)\exp\left[-n\left(\frac{1/n^2-\lambda^2}{2}+2\lambda-\frac{3}{2}\right)\right].
\end{split}
\end{equation}
The splitting into three integrals has the advantage to make simpler estimates for each of the singular terms

 For the first term, the ${\rm arcsinh}$ remains finite and  takes its maximum at $\lambda=1/2$, i.e., we can replace it by ${\rm arcsinh}[2/(2-\sqrt{2})]$ for any $n\geq2$. For the exponential term we look for the maximum of the function
\begin{equation}\label{fdef}
F(\lambda)=\frac{\lambda^2}{2}-2\lambda+\frac{3}{2}+\frac{n-1}{2n}{\rm ln}\left(\lambda^2+\frac{1}{n^2}\right)
\end{equation}
on the interval $[-1/n,1/2]$. Its first derivative 
\begin{equation}
F'(\lambda)=\frac{(\lambda-2)(\lambda^2+1/n^2)+(n-1)\lambda/n}{\lambda^2+1/n^2}
\end{equation}
is strictly increasing for all integer $n$ as can be checked by taking the derivative of its numerator. Since the derivative crosses the real axis once and is negative at the lower bound and positive at the upper limit, the maximum of $F$ is achieved either at $\lambda=-1/n$ or $\lambda=1/2$. The former decreases like $-{\rm ln}(n)$ for very large $n$ while the latter approaches the negative constant $5/8-{\rm ln}(2)\approx-0.07$. Therefore this term is exponentially bounded.

The function~\eqref{fdef} plays also a crucial role in the integration over the other two intervals. Thus, also their maxima are taken at their upper limits. For instance, the exponential function in the third integral over $[1,1+1/(\sqrt{3}n)]$ can be estimated from above at $\lambda=1+1/(\sqrt{3}n)$ like
\begin{equation}
\exp\left[-n\left(\frac{1/n^2-\lambda^2}{2}+2\lambda-\frac{3}{2}\right)\right]\left(\frac{1}{n^2}+\lambda^2\right)^{(n-1)/2}< 24,
\end{equation}
where we have gone slightly further for a simple expression on the right hand side. Thence, we are left by an integral over the ${\rm arcsinh}$ function, the lower limit $\lambda=1$ being the closest point to its singularity $\lambda=\sqrt{1-1/n^2}$. Plugging $\lambda=1$ in, this function behaves like ${\rm ln}(n)$ for large $n$. Yet, the length of the interval is $1/(\sqrt{3}n)$. Therefore, this integral behaves like ${\rm ln}(n)/n$ and converges to zero for $n\to\infty$.

The last integral we have to evaluate is the second part in~\eqref{apprx2.b}. Again the maximum of the exponential function is taken at its upper limit which is this time $\lambda=1$. Since we are in the proximity of the cut of the ${\rm arcsinh}$ function and the saddle point, we cannot easily set any of the terms in the integral to a constant. One term can be done, nevertheless. The factor $\left(1/n^2+\lambda^2\right)^{(n-1)/2}$ is split into $\lambda^{n}$ and  $\left(1+1/(n^2\lambda^2)\right)^{(n-1)/2}/\lambda$. The latter has the upper bound $2\left(1+4/n^2\right)^{(n-1)/2}\leq 2\exp[2(n-1)/n^2]< 16$ for all $n\in\mathbb{N}$. After changing the coordinates to $\lambda=\sqrt{1-\delta t/n^{1/3}}$ with $\delta t\in[0,3n^{1/3}/4]$, leading to a factor $1/n^{1/3}$ in the Jacobian, the exponentially growing term in $n$ becomes
\begin{equation}\label{apprx2.c}
\begin{split}
&\exp\left[n\left(\frac{\lambda^2}{2}-2\lambda+\frac{3}{2}\right)\right]\lambda^{n}\\
=&\exp\left[-n^{2/3}\frac{\delta t}{2}-2n\left(\sqrt{1-\frac{\delta t}{n^{1/3}}}-1\right)+\frac{n}{2}{\rm ln}\left(1-\frac{\delta t}{n^{1/3}}\right)\right].
\end{split}
\end{equation}
Next, we employ the Taylor series of the two terms
\begin{equation}
\begin{split}
2\left(\sqrt{1-\frac{\delta t}{n^{1/3}}}-1\right)=&-\sum_{j=1}^\infty\frac{\Gamma[j-1/2]}{\sqrt{\pi}j!}\left(\frac{\delta t}{n^{1/3}}\right)^j,\\
 \frac{1}{2}{\rm ln}\left(1-\frac{\delta t}{n^{1/3}}\right)=&-\sum_{j=1}^\infty\frac{1}{2j}\left(\frac{\delta t}{n^{1/3}}\right)^j.
\end{split}
\end{equation}
Their difference yields a Taylor series whose coefficients are all negative for $j\geq3$. The coefficient for the second order is equal to zero and for the first order it cancels with the linear term in the exponent~\eqref{apprx2.c}. Hence, the exponentially growing term in $n$ can be bounded by the exponential function $\exp[-\delta t^3/24]$. Because of  the logarithmic behaviour of the ${\rm arcsinh}$ function, the whole contribution of the integral over $\lambda\in[1/2,1]$ is proportional to ${\rm ln}(n)/n^{1/3}$.

Collecting all three contributions of the integral~\eqref{apprx2.b}, we have with three $n$-independent constants $\gamma_1$, $\gamma_2$, and $\gamma_3$ the following estimate
\begin{equation}
\begin{split}
|\mathcal{I}_{\mathcal{C}_2}|\leq& {\rm const.}\biggl(\gamma_1\exp[-0.07 n]+\gamma_2\frac{{\rm ln}(n)}{n^{1/3}}+\gamma_3\frac{{\rm ln}(n)}{n}\biggl).
\end{split}
\end{equation}
Therefore, this part of the integration over the whole contour also vanishes at least like ${\rm ln}(n)/n^{1/3}$, which completes the proof for $|x|=2$.

\subsection{The Case $2<|x|$}

Again, we first consider the simpler case $\tilde{a}_2x>0$ which is now only given by the parametrization $z=t+z_+/|z_-|$ with $t\in\mathbb{R}$. For the function $F$, we obtain
\begin{equation}\label{Fxout}
\begin{split}
\left|F\left(t+\frac{z_+}{|z_-|}\right)\right|=&\exp\left[-n\left(\frac{|z_-|^2t^2-|z_+|^2}{2}+|xz_+|\right)\right]\left(t^2+\frac{|z_+|^2}{|z_-|^2}\right)^{n/2}\\
\leq& \exp\left[-n\frac{1-|z_-|^2}{4|z_+|^2}t^2+n\left(\frac{|z_+|^2}{2}-|xz_+|\right)\right]\left(\frac{|z_+|}{|z_-|}\right)^{n}
\end{split}
\end{equation}
because the maximum of the quotient of the first and the second line is taken at $t=0$. When noticing that $|x|=|z_+|+|z_-|$ and $|z_+z_-|=1$, the ratio of the $z$ and $z'$ integrand is always exponentially smaller than unity, i.e.,
\begin{equation}\label{Fratioout}
\begin{split}
\left|\frac{F\left(t+z_+/|z_-|\right)}{F(e^{i\varphi})}\right|\leq& \exp\left[-n\frac{1-|z_-|^2}{4|z_+|^2}t^2-n\frac{|z_+|^2-|z_-|^2}{2}\right]\left(\frac{|z_+|}{|z_-|}\right)^{n}\\
=&\exp[-n(\sinh[\tilde{\theta}]-\tilde{\theta})]
\end{split}
\end{equation}
with $\tilde{\theta}={\rm ln}(|z_+|/|z_-|)>0$. Additionally, the denominator $1/|e^{i\varphi}-t-z_+/|z_-||$  has the upper bound $1/\sqrt{|z_+|/|z_-|-1}$ as does the function $K\omega^{(n)}\left(\tilde{a}_2(t+\frac{z_+}{|z_-|})\right)$ by $c|z_-|/|\tilde{a}_2z_+|$. Therefore, the whole integrand is uniformly bounded from above by $e^{-n \gamma}$ with $\gamma>0$. Hence, the integral $\mathcal{I}$ vanishes exponentially.

Nothing changes for the situation of $\tilde{a}_2x<0$ when considering the part $\mathcal{C}_1=\{z=t+z_+/|z_-|\,|t\in\mathbb{R}\ {\rm and}\ |t|>1/n\}$ of the whole contour $\mathcal{C}_z=\mathcal{C}_1\cup\mathcal{C}_2\cup\mathcal{C}_3$. Indeed the additional exponential term for the bound of $K\omega^{(n)}$ is inferior to the Gaussian in Eq.~\eqref{Fratioout} and does not even scale with $n$.

We would like to point out that due to a different radius $J\omega^{(n)}$ has now a slightly smaller bound $Ce^{|z_-\tilde{a}_1|}$. Apart from this rescaling everything remains the same for that part.

 The other part $\mathcal{C}_3$ parallel to the real axis with $z=t-i\sign(\tilde{a}_2)/n$ with $t\in[-1/n,1/n]$ leads to the ratio
\begin{equation}
\begin{split}
\left|\frac{F\left(t-i\sign(\tilde{a}_2)/n\right)}{F(e^{i\varphi})}\right|\leq& \exp\left[-n|z_-|^2\left(\frac{t^2-1/n^2}{2}-\frac{|x|}{n|z_-|}\right)+n\frac{2+|z_-|^2}{2}\right]\\
&\times\left(t^2+\frac{1}{n^2}\right)^{n/2}\\
=&\exp\left[n\frac{2+(1-t^2)|z_-|^2}{2}+|z_-x|+\frac{|z_-|}{2n}\right]\left(t^2+\frac{1}{n^2}\right)^{n/2}.
\end{split}
\end{equation}
The function $K\omega^{(n)}\left(\tilde{a}_2(t-i\sign(\tilde{a}_2)/n)\right)$ has again the bound $c/(|\tilde{a}_2|\sqrt{t^2+1/n^2})$, and for  the denominator $1/|e^{i\varphi}-t-i\sign(\tilde{a}_2)/n|$ it is given by $\sqrt{2}/\sqrt{1-2/n^2}\leq 2$. Collecting everything we have
\begin{equation}
\begin{split}
|\mathcal{I}_{\mathcal{C}_3}|\leq &{\rm const.}\int_{-1/n}^{1/n}\frac{dt}{2\pi}\exp\left[n\frac{2+(1-t^2)|z_-|^2}{2}+|z_-x|+\frac{|z_-|}{2n}\right]\left(t^2+\frac{1}{n^2}\right)^{(n-1)/2}\\
\leq&{\rm const.}\ \exp\left[n\frac{2+|z_-|^2}{2}+|z_-x|+\frac{|z_-|}{2n}\right]\left(\frac{2}{n^2}\right)^{n/2}.
\end{split}
\end{equation}
When choosing $n>\sqrt{2}\exp[1+|z_-|^2/2]$ we see that this contribution is also exponentially suppressed in the limit $n\to\infty$ as it has been the case for $|x|<2$.

Eventually, we consider the parts in $\mathcal{C}_2$ that are parallel to the imaginary axis, especially we consider the parametrization $z=\sign(t)/n-i \sign(x)\lambda$ with $\lambda=|tz_+/z_-|-(1-|t|)/n\in[-1/n,|z_+/z_-|]$ and, thus, $t\in[-1,1]$. This means that $K\omega^{(n)}\left(\tilde{a}_2(\sign(t)/n-i\sign(x)\lambda)\right)$ is replaced by $c/(|\tilde{a}_2|\sqrt{\lambda^2+1/n^2})$ and $1/|e^{i\varphi}-\sign(t)/n+i \sign(x)\lambda|$ by $\sqrt{2}/\sqrt{(\varphi-\vartheta(t))^2/\pi^2+(1-\sqrt{\lambda^2+1/n^2})^2}$ with $\sign(t)/n-i \sign(x)\lambda=\sqrt{\lambda^2+1/n^2}e^{\vartheta(t)}$. The additional exponential term $(|z|+1)\exp[\alpha|z|]$ for $K\omega^{(n)}$ can be again absorbed into the constant $c$ due to the compact support of this part of the contour, in particular the maximum can be estimated from above by the $n$-independent constant $(|\tilde{a}_2|\sqrt{1+|z_+/z_-|^2}+1)\exp[\alpha|\tilde{a}_2|\sqrt{1+|z_+/z_-|^2}]$.

The angle $\vartheta(t)$ obviously depends on the sign of $t$, too. However, this has no effect since for each of the two components of $\mathcal{C}_2$, we need to shift $\varphi$ such that $\varphi-\vartheta(\lambda)\in]-\pi,\pi]$ to apply the aforementioned approximation, cf., Eq.~\eqref{estimate.ratio}. Thence, both parts yield the same contribution after employing the approximations and integrating over $\varphi$. Then, we arrive at
\begin{equation}
\begin{split}
|\mathcal{I}_{\mathcal{C}_2}|\leq& {\rm const.}\int_{-1/n}^{|z_+/z_-|}\frac{d\lambda\ \left(\lambda^2+\frac{1}{n^2}\right)^{(n-1)/2}}{2\pi(|z_+/z_-|+1/n)}{\rm arcsinh}\left[\frac{1}{|1-\sqrt{\lambda^2+1/n^2}|}\right]\\
&\times\exp\left[-n\left(|z_-|^2\frac{1/n^2-\lambda^2}{2}+(1+|z_-|^2)\lambda\right)+n\frac{2+|z_-|^2}{2}\right],
\end{split}
\end{equation}
where we have also exploited
\begin{equation}
\begin{split}
\left|F\left(\sign(t)t_0-i\sign(x) \lambda\right)\right|=&\exp\left[-n|z_-|^2\left(\frac{1/n^2-\lambda^2}{2}+\frac{|x|}{|z_-|}\lambda\right)\right]\left(\frac{1}{n^2}+\lambda^2\right)^{n/2}
\end{split}
\end{equation}
and Eq.~\eqref{zprime.est}.

The exponentially growing term in $n$ that depends on $\lambda$ takes two extrema at $\lambda=|z_+|/|z_-|$ and $\lambda=1$. We, here, neglected $1/n^2$ corrections from the $(\lambda^2+1/n^2)^{n/2}$ term. The first extrema is a minimum while the second is a maximum. Due to that, we shift $\lambda\to 1+\lambda$ leading to
\begin{equation}
\begin{split}
|\mathcal{I}_{\mathcal{C}_2}|\leq& {\rm const.}\int_{-1/n-1}^{|z_+/z_-|-1}\frac{d\lambda}{2\pi}{\rm arcsinh}\left[\frac{1}{|1-\sqrt{(1+\lambda)^2+1/n^2}|}\right]\\
&\times\left((\lambda+1)^2+\frac{1}{n^2}\right)^{(n-1)/2}\exp\left[n\left(|z_-|^2\frac{\lambda^2}{2}-\lambda\right)\right]\\
\leq&{\rm const.}\biggl(\int_{-1/n-1}^{-1/2}+\int_{-1/2}^{(1-|z_-|^2)/(1+|z_-|^2)}+\int_{(1-|z_-|^2)/(1+|z_-|^2)}^{|z_+/z_-|-1}\biggl)\frac{d\lambda}{2\pi}\\
&\times{\rm arcsinh}\left[\frac{1}{|1-\sqrt{(1+\lambda)^2+1/n^2}|}\right]\left((\lambda+1)^2+\frac{1}{n^2}\right)^{(n-1)/2}\\
&\times\exp\left[n\left(|z_-|^2\frac{\lambda^2}{2}-\lambda\right)\right],
\end{split}
\end{equation}
where we have additionally dropped $\exp[-|z_-|^2/2n]\leq1$. Since we stay away from $\lambda=0$ in the integration over $[-1/n-1,-1/2]$, the ${\rm arcsinh}$ is bounded by ${\rm arcsinh}[2/(2-\sqrt{2})]$ for any $n\geq2$ and the  exponentially growing term in $n$ takes its maximum at $\lambda=-1/2$, compare this with the case $|x|=2$ since the discussion is exactly the same.

Similarly for $[(1-|z_-|^2)/(1+|z_-|^2),|z_+/z_-|-1]$ both terms are bounded, this time by their lower boundary at $\lambda=(1-|z_-|^2)/(1+|z_-|^2)$.
Due to this we arrive at
\begin{equation}\label{pxg2.a}
\begin{split}
|\mathcal{I}_{\mathcal{C}_2}|\leq&{\rm const.}\biggl(n\gamma_1\exp\left[\frac{n}{2}\left(\frac{|z_-|^2}{4}-1\right)\right]+\gamma_2\left(\frac{4}{(|z_-|^2+1)^2}+\frac{1}{n^2}\right)^{n/2}\\
&\exp\left[\frac{n}{2}\frac{|z_-|^6+|z_-|^2-2}{(1+|z_-|^2)^2}\right]+\int_{-\sqrt{n}/2}^{\sqrt{n}(1-|z_-|^2)/(1+|z_-|^2)}\frac{d\lambda}{2\pi\sqrt{n}}e^{|z_-|^2\frac{\lambda^2}{2}-\sqrt{n}\lambda}\\
&\times{\rm arcsinh}\left[\frac{1}{|1-\sqrt{(1+\lambda/\sqrt{n})^2+1/n^2}|}\right]\left(\left(\frac{\lambda}{\sqrt{n}}+1\right)^2+\frac{1}{n^2}\right)^{(n-1)/2}
\end{split}
\end{equation}
with two constants $\gamma_1,\gamma_2>0$ and $n$ large enough.
The first two terms vanish exponentially with $n\to\infty$ because $|z_-|<1$. For the second term, it is worthwhile to note that between the two aforementioned extrema the exponential is strictly decreasing and the maximum at $\lambda=(1-|z_-|^2)/(1+|z_-|^2)$ lies inside the integration domain so that indeed its behaviour is $e^{-n\gamma}$ with $\gamma>0$.

In the third term of~\eqref{pxg2.a}, we have rescaled $\lambda\to\lambda/\sqrt{n}$. Therein, the function ${\rm arcsinh}|x|$ can be replaced by its upper bound ${\rm arcsinh}\left[1\right]+|{\rm ln}|x|\ |$, and the factor of $1/\sqrt{(\lambda/\sqrt{n}+1)^2+1/n^2}$ is bounded by $2$,
\begin{equation}
\begin{split}
|\mathcal{I}_{\mathcal{C}_2}|\leq&{\rm const.}\biggl(n\gamma_1\exp\left[\frac{n}{2}\left(\frac{|z_-|^2}{4}-1\right)\right]+\gamma_2\left(\frac{4}{(|z_-|^2+1)^2}+\frac{1}{n^2}\right)^{n/2}\\
&\exp\left[\frac{n}{2}\frac{|z_-|^6+|z_-|^2-2}{(1+|z_-|^2)^2}\right]+\int_{-\sqrt{n}/2}^{\sqrt{n}(1-|z_-|^2)/(1+|z_-|^2)}\frac{d\lambda}{\pi\sqrt{n}}e^{|z_-|^2\frac{\lambda^2}{2}-\sqrt{n}\lambda}\\
&\hspace*{-0.5cm}\times\left({\rm arcsinh}[1]+\left|{\rm ln}\left|1-\sqrt{(1+\lambda/\sqrt{n})^2+1/n^2}\right|\ \right|\right)\left(\left(\frac{\lambda}{\sqrt{n}}+1\right)^2+\frac{1}{n^2}\right)^{n/2}\biggl),
\end{split}
\end{equation}
The behaviour of the exponentially growing part in $n$ is better understood with the aid of the function
\begin{equation}
\begin{split}
L\left(\frac{\lambda}{\sqrt{n}}\right)=\frac{|z_-|^2+1}{4}\frac{\lambda^2}{n}-\frac{\lambda}{\sqrt{n}}+{\rm ln}\left(1+\frac{\lambda}{\sqrt{n}}\right)
\end{split}
\end{equation}
which is monotonously growing in $n>0$ for all fixed $\lambda\in]-\sqrt{n}/2,\sqrt{n}(1-|z_-|)/(1+|z_-|)[$. One can show this via its first derivative in $n$ leading to $-\lambda/(2n^{3/2})L'(\lambda/\sqrt{n})$ and noticing that $L'(\lambda/\sqrt{n})$ is positive/negative for negative/positive $\lambda$ in the allowed interval. Moreover, $L(0)=0$ so that we have
\begin{equation}
\begin{split}
\frac{|z_-|^2}{2}\frac{\lambda^2}{n}-\frac{\lambda}{\sqrt{n}}+{\rm ln}\left(1+\frac{\lambda}{\sqrt{n}}\right)\leq \frac{|z_-|^2-1}{4}\frac{\lambda^2}{n}\leq 0.
\end{split}
\end{equation}
With the aid of this identity and by noticing that $(1+1/(\sqrt{n}\lambda+n)^2)^{n/2}\leq (1+4/n^2)^{n/2}\leq e^{2}$, we have
\begin{equation}
\begin{split}
|\mathcal{I}_{\mathcal{C}_2}|\leq&{\rm const.}\biggl(n\gamma_1\exp\left[\frac{n}{2}\left(\frac{|z_-|^2}{4}-1\right)\right]+\gamma_2\left(\frac{4}{(|z_-|^2+1)^2}+\frac{1}{n^2}\right)^{n/2}\\
&\exp\left[\frac{n}{2}\frac{|z_-|^6+|z_-|^2-2}{(1+|z_-|^2)^2}\right]+\int_{-\sqrt{n}/2}^{\sqrt{n}(1-|z_-|^2)/(1+|z_-|^2)}\frac{d\lambda}{\pi\sqrt{n}}\\
&\hspace*{-0.5cm}\times\left({\rm arcsinh}[1]+\left|{\rm ln}\left|1-\sqrt{(1+\lambda/\sqrt{n})^2+1/n^2}\right|\ \right|\right)\exp\left[\frac{|z_-|^2-1}{4}\lambda^2+2\right]\biggl).
\end{split}
\end{equation}
The logarithm  becomes infinitely large when the argument becomes zero, namely at $\lambda=\sqrt{n-1/n}-\sqrt{n}$ which lies close to $0$ for large $n$. Hence, performing an additional shift $\lambda\to\lambda+\sqrt{n}(\sqrt{1-1/n^2}-1)$ exhibits the divergent behaviour of the logarithm as $|{\rm ln}|\lambda/\sqrt{n}|\ |$ which grows logarithmically with $n$ but is still integrable. We underline that the shift is for large $n$ only of order $n^{-3/2}$ which does not destroy our convergence generating Gaussian in the integral.

When summarizing the above discussion, the integral $\mathcal{I}$ approaches zero at least like ${\rm ln}(n)/\sqrt{n}$. This closes the proof of Theorem~\ref{thm:main}.

\section{Conclusions and a Conjecture}\label{sec:conclusio}

We have proven the hard edge statistics of the singular values of the complex square matrix  $G$ drawn from a P\'olya ensemble, see Proposition~\ref{prop:HEKPolya}, and of the eigenvalues of the product matrix $G(H-nx\eins_n)G^*$, where $x\in\mathbb{R}$ is fixed and $H$ is drawn from a GUE. The conditions for the P\'olya ensemble are not extremely strong, cf., Eq.~\eqref{assump1} and~\eqref{assump2}, and seem to be intimately related to the fact that $G$ has a hard edge limit at the origin. We have not proven that they correspond exactly to the existence of a hard edge for the singular values of $G$ at the origin but could see that indeed for ensembles that have a spectral gap about the origin those conditions are not satisfied. For example, the inverse Ginibre matrix is such a random matrix ensemble. It would be interesting whether one can completely reduce the condition to the fact of the existence of a hard edge limit at the origin.

When considering our main result Theorem~\ref{thm:main}, one readily notices two facts. Firstly, the hard edge statistics of $G$ plays a bigger role than the local spectral statistics of $H-nx\eins_n$. In fact, we obtain the hard edge kernel of $G$ much earlier than expected. Certainly, the matrix $H$ plays a subleading role when $|x|\gg1$ since the product matrix is dominated by $-n x GG^*$ which leads to nothing else than the rescaled singular value statistics of $G$. Yet, the kernel for $G(H-nx\eins_n)G^*$ is already identical to the one of $GG^*$, apart from a rescaling, when $x$ lies outside the support of the spectral density of $H$, which is the Wigner semi-circle. Thus, we do not need at all that $|x|$ has to be immensely large.

Another remarkable point is that the local soft edge statistics of the GUE, which is the Airy kernel~\cite{Forresterbook}, becomes only visible due to a change of the rate of convergence. Although we have not looked for the optimal rates, we have seen in our proof that it drastically changes when $x$ is on one of the two soft edges of the GUE. This has been corroborated by numerical simulations where $G$ has been drawn from a Ginibre ensemble.

Interestingly, the Stieltjes transform (or Green function) of the GUE plays a much more important role and determines the rescaling of two terms in the kernel of the product matrix. In principle one can do the same calculation when replacing $H$ by any other additive P\'olya ensemble like a Laguerre ensemble~\cite{Forresterbook} or a Muttalib--Borodin ensemble~\cite{Borodin1998} (where $\omega(\lambda)=e^{-\nu \lambda}\exp[-e^{-\alpha\lambda}]$ with $\alpha,\nu>0$) or even no additive P\'olya ensemble at all like a shifted Jacobi (truncated unitary matrices)~\cite{KKS2016} or Cauchy-Lorentz~\cite{KK2016,KK2019} ensemble. For instance for the aforementioned  additive P\'olya ensembles, one can show the same formula with some modifications of the proof, where the Green function of the GUE has to be replaced by the one of the considered ensemble. Due to this, we believe that the result is much more general than stated above. Therefore, we would like to conclude with the following conjecture. Its corroboration with its proof for the considered random matrix $G(H-nx\eins_n)G^*$ with $H$ drawn from a GUE has been one of our major motivations of writing the present work.

\begin{conjecture}[Hard Edge Statistics of Products involving P\'olya Ensembles]\label{conjecture}\

Let $G$ be a complex square matrix satisfying the conditions~\eqref{assump1} and~\eqref{assump2} and $H$ is a random Hermitian matrix independently drawn from a polynomial ensemble with a limiting Green function $G(x)$ that is finite at $x=0$, i.e., $0<|G(0)|<\infty$. Then, the limiting hard edge kernel of the eigenvalues for the product matrix $GHG^*$ is equal to
\begin{equation}
\begin{split}
K_\infty(\tilde{a}_1,\tilde{a}_2)=&\Theta(-{\rm Re}[ G(0)]\tilde{a}_2)\,|{\rm Re}[ G(0)]|\,K_\infty^{(G)}( \sign(\tilde{a}_2)|{\rm Re}[ G(0)]|\tilde{a}_1,|{\rm Re}[ G(0)]\tilde{a}_2|)\\
&+{\rm Im}[G(0)]\int_{-1}^{1}\frac{dt}{2\pi} J\omega^{(\infty)}\left(-i\tilde{a}_1({\rm Im}[G(0)]t-i{\rm Re}[ G(0)])\right)\\
&\times K\omega^{(\infty)}\left(\tilde{a}_2({\rm Im}[G(0)]t-i{\rm Re}[ G(0)])\right),
\end{split}
\end{equation}
where we define $G(0)=\lim_{\epsilon\searrow 0}G(i\epsilon)$ such that ${\rm Im}[G(0)]\geq0$.
\end{conjecture}

The situations of $|G(0)|=0$ and $|G(0)|=\infty$ need to be excluded. Both indicate a different scaling behaviour; the first one  because the saddle points come very close to the cut of $K\omega^{(\infty)}$ and the second one since the exponential bounds of $J\omega^{(\infty)}$ and $K\omega^{(\infty)}$ become important. The case $|G(0)|=0$ is exactly the case what happens in~\cite[Theorem 3 (ii) and (iii)]{Liu2017} which explains why another kernel has been found. However, the case~\cite[Theorem 3 (i)]{Liu2017} again agrees with our findings since the Green function of the Hermitian matrix $H$ acquires an imaginary part, in particular the macroscopic level density is non-zero.

We have not proven this form of the theorem since it is certainly technically  more involved. It might be, as already mentioned, that the condition for $G$ can be relaxed to an existing hard edge kernel and that $H_0$ can be any Hermitian matrix, even fixed. However, then the limiting  determinantal point process for the hard edge is not obvious, which is guaranteed for a polynomial ensemble at finite matrix dimension, see~\cite{KK2016,KK2019}.

\section*{Acknowledgements}

I am grateful for the fruitful discussions with  Dang-Zheng Liu and Holger K\"osters who also read the first draft and gave their feedback. Moreover, I would like to thank the anonymous referee for their helpful comments.

Additionally, I would like dedicate the present article to my wife for wonderful ten years of marriage.

\end{document}